  \documentclass[12pt]{amsart}
 \usepackage{txfonts}
 \usepackage{mathrsfs}
 \usepackage{amsfonts}

 \usepackage[dvips]{graphicx}

 \usepackage{epsfig, graphicx, graphics, color, amsmath,amssymb, amsfonts, amstext,amsthm, latexsym}

   \newtheorem{lemma}{Lemma}[section]

   \newtheorem{prop}[lemma]{Proposition}
   
   \newtheorem{definition}[lemma]{Definition}
   \newtheorem{example}[lemma]{Example}

\newcommand{\s}{\sigma}

\newcommand{\e}{\epsilon}
\renewcommand{\a}{\alpha}
\renewcommand{\b}{\beta}

\newcommand{\D}{\Delta}
\newcommand{\p}{\partial}

\renewcommand{\phi}{\varphi}
\newcommand{\N}{{\mathbb N}}
\newcommand{\R}{{\mathbb R}}

\newcommand{\F}{{\mathcal{F}}}
\newcommand{\B}{{\mathcal{B}}}

\makeatletter
    \newcommand\figcaption{\def\@captype{figure}\caption}
    \newcommand\tabcaption{\def\@captype{table}\caption}
\makeatother

\title[Quantifying Uncertainties in Complex Systems]
{Quantifying Model Uncertainties \\in Complex  Systems}

\author{Jiarui Yang and Jinqiao Duan  }
\address[Jiarui Yang and Jinqiao Duan ]
{Department of Applied Mathematics\\
 Illinois Institute of Technology\\
   Chicago, IL 60616, USA }
\email{J. Yang: jyang31@iit.edu; J. Duan: duan@iit.edu}

 \date{\today}

\keywords{Model uncertainty, parameter estimation,  Brownian motion (BM), fractional
Brownian motion (fBM),
 L\'evy motion (LM), Hurst parameter, characteristic exponent, stochastic differential equations (SDEs) \\
This work was partially supported by NSF grants 0620539 and 0731201, and by an open research grant of the State Key Laboratory for Nonlinear Mechanics, China. }

\begin{document}

\def \a {\alpha}
\def \b {\beta}
\def \m {\mu}
\def \O {\Omega}
\def \e {\acute{e}}
\def \D {\Delta}
\def \f {\forall}
\def \l {\lambda}
\def \B {\mathscr{B}}
\def \R {{\mathbb{R}}}
\def \F {\mathcal {F}}
\def \C {{\bf C}}
\def \E {{\bf E}}
\def \N {{\mathbb{N}}}
\def \P {\mathbb{P}}
\def \Q {{\bf Q}}
\def \p{\partial}
\def \t{\theta}
\def \s {\sigma}
\def \ba{\begin{eqnarray*}}
\def \ea{\end{eqnarray*}}
\def \bc{\begin{center}}
\def \ec{\end{center}}
\def \bl{\begin{flushleft}}
\def \el{\end{flushleft}}
\def \br{\begin{flushright}}
\def \er{\end{flushright}}
\def \bd{\begin{document}}
\def \ed{\end{document}}

\begin{abstract}
Uncertainties are abundant in complex systems. Appropriate mathematical models
for these systems thus contain random effects or noises. The models
are often in the form of stochastic differential equations, with
some parameters to be determined by observations. The stochastic
differential equations may be driven by Brownian motion, fractional
Brownian motion, or L\'evy motion.

After a brief overview of recent advances in estimating parameters in stochastic differential equations, various numerical algorithms for computing parameters are implemented.
 The numerical simulation results
are shown to be   consistent with   theoretical analysis.
Moreover, for fractional Brownian
motion and $\alpha-$stable L\'evy motion, several algorithms  are
reviewed and implemented to numerically estimate the Hurst parameter $H$ and characteristic
exponent $\alpha$.

\end{abstract}

\maketitle

\section{Introduction}
\label{intro}

Since random fluctuations   are common   in the
real world, mathematical models for complex systems are often subject to uncertainties,
such as fluctuating forces, uncertain parameters, or random boundary conditions \cite{Moss, Horst, Gar, VanKampen3, WaymireDuan, Wong}.  Uncertainties may also be caused by the lack of knowledge of some
    physical, chemical or biological mechanisms  that are   not
    well understood,  and thus are not appropriately
    represented (or missed completely) in the mathematical models \cite{ChenDuan, Kantz, Palmer2, Wilks, Williams}.

\medskip

Although these fluctuations and  unrepresented     mechanisms  may
be very small or very fast, their long-term impact on the system
evolution
  may be delicate \cite{Arnold, Horst, Gar}   or even profound. This kind of
delicate impact on the overall evolution of dynamical  systems has
been observed in, for example, stochastic bifurcation
\cite{Crauel, CarLanRob01, Horst}, stochastic resonance \cite{Imkeller},
 and  noise-induced pattern formation \cite{Gar, BlomkerStani}.
Thus taking stochastic effects   into account is of
central importance for   mathematical modeling of
complex systems under uncertainty. Stochastic     differential equations (SDEs) are appropriate models for many of these systems \cite{Arnold, DaPrato, Roz, WaymireDuan}.

\medskip

For example,   the Langevin type models   are stochastic
differential equations describing various phenomena in physics,
biology,  and  other   fields. SDEs are used to model various price
processes, exchange rates, and interest rates, among others, in
finance. Noises in these SDEs may be modeled as a generalized time
derivative of some distinguished stochastic processes, such as
Brownian motion (BM), L\'evy motion (LM) or fractional Brownian
motion (fBM) \cite{DuanLiWang}. Usually we choose   different noise
processes according to the statistical property of the observational
data. For example, if the data has long-range dependence, we
consider fractional Brownian motion rather than Brownian motion. If
the data has considerable discrepancy with Gaussianity or normality,
L\'evy motion may  be an appropriate choice. In building these SDE
models, some parameters appear, as we do not know certain quantities
exactly.

Based on the choice of noise processes, different mathematical
techniques are needed in       estimating the parameters in SDEs
with Brownian motion, L\'evy motion, or fractional Brownain motion.

\medskip

In this article, we are interested in estimating and computing
  parameters contained in stochastic differential
equations, so that we obtain computational models useful for
investigating complex dynamics under uncertainty. We first review
recent theoretical results in estimating parameters in SDEs,
including statistical properties and   convergence of various
estimates. Then we develop and implement numerical algorithms in
approximating these parameters.

\medskip

Theoretical results on parameter estimations for SDEs driven by
Brownian motion are relatively well developed
(\cite{Alizadeh2002, Davis2000,Genon-Catalot1999, Pearson1994}), and
various numerical simulations for these parameter estimates
(\cite{Ait2002, Ait-Sahalia2003, Pearson1994, Nicolau2004})  are
implemented. So, in Section 2 below, we do not present such
numerical results. Instead, we will concentrate on numerical
algorithms for parameter estimations in SDEs driven by fractional
Brownian motion and L\'evy motion in Section 3 and 4, respectively.

\bigskip

This paper is organized as follows. In
\textbf{Section 2}, we consider parameter estimation for SDEs with
  Brownian motion $B_t$. We   present a brief overview of
some available techniques on estimating parameters in these
stochastic differential equations with continuous-time or
discrete-time observations. In fact, we present  results about how
to estimate parameters in diffusion terms and  drift terms, given
continuous observations and discrete observations, respectively.

\medskip

In \textbf{Section 3}, we consider parameter estimation for  SDEs
driven by fractional Brownian motion $B_t^H$ with Hurst parameter $H$.
After discussing basic properties of fBM, we consider
parameter estimation methods for Hurst parameter $H$ from given fBM
data. Then, we compare the convergence rate of each method by
comparing estimates computed with hypothetic data. Unlike the case
of SDEs with Brownian motion,
    there is no general formula for the estimate of the
parameter in the drift (or diffusion) coefficient of a stochastic
differential equation driven by fBM.   We discuss different
estimates associated with different models and discuss the
statistical properties respectively. We develop and implement
numerical simulation methods for these estimates.

\medskip

Finally, in \textbf{Section 4}, for   stochastic differential
equation with (non-Gaussian)  $\alpha-$stable L\'evy motion
$L_t^\alpha$, we consider estimates and their numerical
implementation for  parameter $\alpha$ and other parameters in the
drift or diffusion coefficients.



\section{Quantifying Uncertainties in SDEs Driven by Brownian motion}

In this section, we   consider a   scalar diffusion
process $X_t\in \mathbb{R}^d, 0\leq t \leq T$ satisfying the
following stochastic differential equation
\begin{equation}\label{generalequationBM}
dX_t=\mu(\theta, t, X_t)dt+\sigma(\vartheta, t, X_t)d B_t, ~~~
X_0=\zeta
\end{equation}
where $B_t$ is a m-dimensional Brownian motion, $\theta\in \Theta$ a
compact subset of $\mathbb{R}^p$ and $\vartheta\in \Xi$ a compact
subset of $\mathbb{R}^q$ are unknown parameters which are to be
estimated on the basis of observations. Here $\mu: \Theta\times [0,
T]\times \mathbb{R}^d\rightarrow \mathbb{R}^d$,   the drift
coefficient, and $\sigma: \Xi\times [0, T]\times
\mathbb{R}^d\rightarrow \mathbb{R}^{d\times m}$,   the diffusion coefficient,
are usually known functions but with unknown parameters $\theta$ and $\vartheta$.\\

Some remarks are in order here.
 \begin{itemize}
            \item Under local Lipschitz and the sub-linear growth conditions on the
coefficients $\mu$ and $\sigma$, there exists a unique strong
solution of the above stochastic differential equation (see
\cite{Kutoyants1084b} or \cite{Mao1995}) and this is an universal
assumption for all results we discuss below.
            \item The diffusion coefficient $\sigma$ is almost surely determined by the
process, i.e., it can be estimated without any error if observed
continuously throughout a time interval (see ~\cite{Genon1994} and
\cite{Doob1953}).
            \item The diffusion matrix defined by $\Sigma(\vartheta, t, X_t)\equiv\sigma(\vartheta, t, X_t)\sigma(\vartheta, t, X_t)^T$
            plays an important role on parameter
            estimation problems.
          \end{itemize}


\subsection{How to Estimate Parameters Given Continuous Observation}
Since it is not easy to estimate parameters $\theta$ and $\vartheta$
at the same time, usually we simplify our model by assuming one of
those parameters is known and then estimate the other. Moreover,
instead of representing the results of all types of diffusion
processes, we choose to present the conclusion of the most general
one, such as, we prefer the nonhomogeneous case rather than the
homogeneous one or the nonlinear one rather than the linear one.

\subsubsection{Parameter Estimation of Diffusion Terms with Continuous
Observation }

We assume that the unknown parameter $\theta$ in the drift
coefficient is known. Then our model can be simplified as
\begin{equation}\label{generalBMdiffusionCon}
dX_t=\mu(t, X)dt+\sigma(\vartheta, t, X_t)d B_t, ~~~ X_0=\zeta
\end{equation}


Remarks: \begin{itemize}
          \item Different with the model \eqref{generalequationBM}, the drift coefficient $\mu(t, X)$ in model \eqref{generalBMdiffusionCon} is possibly unknown
and maybe related to the whole past of process $X$ instead of $X_t$.
In this case, our model can be easily extended to the non-Markovian
case which is more general than case (1).
          \item If $\mu$ is depending on the unknown parameter $\vartheta$ in model \eqref{generalBMdiffusionCon}, we can also prove the local asymptotic mixed normality property for the
maximum likelihood estimate (MLE) when $\mu(t, X)=\mu(\vartheta,
X_t)$ and $\sigma(\vartheta, t, X_t)=\sigma(\vartheta, X_t)$ (see
\cite{Dohnal1987}).\\
        \end{itemize}

If the diffusion matrix $\Sigma(\vartheta, t, X_t)$ is invertible,
then define a family of contrasts by
\begin{equation}\label{Contrast Function}
U^n(\vartheta)=\frac{1}{n}\sum^n_{i=1}[\log \det \Sigma(\vartheta,
t_{i-1}^n, X_{t_{i-1}^n})+(X^n_i)^T\Sigma(\vartheta,t_{i-1}^n,
X_{t^n_{i-1}})^{-1}X^n_i],
\end{equation}
where
\begin{equation*}
X^n_i=\frac{1}{\delta^n_i}(X_{t^n_i}-X_{t^n_{i-1}}),~\delta_i^n=t^n_i-t^n_{i-1}~,~for~1\leq
i \leq n,
\end{equation*}
and $t^n_i$ is an appropriate partition of [0,T]. However, this
assumption does not always hold. So, we consider a more general
class of contrasts of the form
\begin{equation}\label{General Contrast Function}
U^n(\vartheta)=\frac{1}{n}\sum^n_{i=1}f(\Sigma(\vartheta,
t^n_{i-1},X_{t^n_{i-1}}),X^n_i),
\end{equation}
where $f$ should satisfy certain conditions to obtain the asymptotic
property and consistency property for the estimate generated by
these contrasts below (see \cite{Valentine and
Jacod 1993}). \\

Let $\widehat{\vartheta}_n$ be a minimum contrast estimate
associated with $U^n$, i.e. $\widehat{\vartheta}_n$ satisfies the
following equation
\begin{equation*}
U^n(\widehat{\vartheta}_n)=\min_{\vartheta\in \Xi}U^n(\vartheta).
\end{equation*}
Under some smoothness assumptions on the coefficient $\mu$ and
$\theta$, empirical sampling measure assumption on the sample times
$t^n_i$, and identifiability assumption on the law of the solution
of \eqref{generalBMdiffusionCon},   Genon-Catalot and  Jacod
\cite{Genon1994} have proved that the estimate
$\widehat{\vartheta}_n$ has a local asymptotic mixed normality
property, i.e., $\sqrt{n}(\widehat{\vartheta}_n-\vartheta_0)$ where
$\vartheta_0$ is the true value of the parameter converges in law to N(0, S).\\

Remarks:\begin{itemize}
          \item We do not include the drift coefficient $\mu$ in the contrast $U^n(\vartheta)$ because it
          is possibly unknown. Even if it is known, we still do not
          want it involved since it is a function of the whole past of X and thus is not observable.
          \item If the diffusion matrix $\Sigma$ is invertible, it can be proven that the contrast of form
          \eqref{Contrast Function} is optimal in the class of
          contrasts of type \eqref{General Contrast Function}.
        \end{itemize}

\subsubsection{Parameter Estimation of Drift Terms with Continuous
Observations.}

We assume that the unknown parameter $\vartheta$ in the diffusion
coefficient is known. Then the model \eqref{generalequationBM} can
be simplified as
\begin{equation}\label{generalBMdriftCon}
dX_t=\mu(\theta,t, X_t)dt+\sigma(t, X_t)d B_t, ~~~ X_0=\zeta.
\end{equation}
Since no good result for above general model exists, we introduce
the result for the following
nonhomogeneous diffusion process instead.\\

Consider a real valued diffusion process $\{X_t, t\geq 0\}$
satisfying the following stochastic differential equation:
\begin{equation}\label{generalBMdriftCon2}
dX_t=\mu(\theta, t, X_t)dt+d B_t, ~~~ X_0=\zeta,
\end{equation}
where the drift coefficient function $\mu$ is assumed to be
nonanticipative. Denote the observation of the process by
$X^T_0:=\{X_t, 0\leq t\leq T\}$ and let $P^T_{\theta}$ be the
measure generated by the process $X^T_0$. Then the Radon-Nicodym
derivative (likelihood function) of $P^T_{\theta}$ with respect to
$P^T_{\theta_0}$ where $\theta_0$ is the true value of the parameter
$\theta$ is given by (see \cite{Lipster1977})
\begin{eqnarray*}
&&L_T(\theta):=(dP^T_{\theta}/dP^T_{\theta_0})(X^T_0)\\
&&=\exp\{\int^T_0[\mu(\theta,
t,X_t)-\mu(\theta_0,t,X_t)]dX_t-\frac{1}{2}\int^T_0[\mu^2(\theta,
t,X_t)-\mu^2(\theta_0,t,X_t)]dt\}.
\end{eqnarray*}
So we can get the Maximal Likelihood Estimate (MLE) defined by
\begin{equation*}
\widehat{\theta}_T:=argsup_{\theta\in \Theta}L_T(\theta).
\end{equation*}
Then we can show that the MLE is strongly consistent, i.e.,
$\widehat{\theta}_T\rightarrow \theta_0 ~
P_{\theta_0}-a.s.~as~T\rightarrow \infty$, and converge to a normal
distribution (see Chapter 4 in \cite{Jaya P. N. Bishwal} for more
details).

Remarks:\begin{itemize}
         \item In \cite{Jaya P. N. Bishwal}, Bishwal also proves that
         the MLE and a regular class of Bayes estimates (BE) are
         asymptotically equivalent.
         \item By applying an increasing
transformation as described in \cite{Ait2002},
\begin{equation}
Y_t=g(X)\equiv \int^X\frac{du}{\sigma(u)},
\end{equation}
we can transform the diffusion process $X_t$ defined by
\begin{equation*}
dX_t=\mu(\theta, X_t)dt+\sigma(X_t)dB_t
\end{equation*}
into another diffusion process $Y_t$ defined by
\begin{equation*}
d\tilde{Y_t}=\tilde{\mu}(\theta, \tilde{Y_t})dt +dB_t,
\end{equation*}
where
\begin{equation}
\tilde{\mu}(\theta, y)\equiv \frac{\mu(g^{-1}(y),
\theta)}{\sigma(g^{-1}(y))}-\frac{1}{2}\frac{\partial\sigma
(g^{-1}(y))}{\partial y}.
\end{equation}
Then, we can get the MLE of process $X_t$ by calculating the MLE of
process $Y_t$ according to what we learned in this section (see
\cite{Ait2002} or \cite{Mykland2004} for more details).
       \end{itemize}

\subsection{How to Estimate Parameters given Discrete Observation}
 Given the practical difficulty in obtaining a complete
continuous observation, we now discuss    parameter estimations
with discrete observation.
\subsubsection{Parameter Estimation of Drift Terms with Discrete
Time} In this section, we assume that the unknown parameter
$\vartheta$ in the diffusion coefficient $\sigma$ is known. Then the
model \eqref{generalequationBM} can be simplified as
\begin{equation}\label{generalBMdriftCon}
dX_t=\mu(\theta,t, X_t)dt+\sigma(t, X_t)d B_t, ~~~ X_0=\zeta.
\end{equation}
Ideally, when the transition densities $p(s,x,t,y;\theta)$ of X are
known, we can use the log likelihood function
\begin{equation*}
l_n(\theta)=\sum^n_{i=1} \log p(t_{i-1}, X_{t_{i-1}}, t_i,
X_{t_{i}}; \theta),
\end{equation*}
to compute the LME $\widehat{\theta}$ which is strongly consistent
and asymptotically normally distributed. (see \cite{Billingsley1961}
and \cite{Dacumha1986}, \cite{Le Breton 1976}
and \cite{Robinson1977}).\\

If the transition densities of X are unknown, instead of computing
the log likelihood function $l_n(\theta)$, we would like to use
approximate log-likelihood function which, under some regularity
conditions (see \cite{Hutton1986}), is given by
\begin{equation*}
l_T(\theta)=\int^T_0\frac{\mu(\theta,t,X_t)}{\sigma^2(t,X_t)}dX_t-\frac{1}{2}\int^T_0\frac{\mu^2(\theta,t,X_t)}{\sigma^2(t,X_t)}dt
\end{equation*}
to approximate the log-likelihood function based on continuous
observations (see \cite{prakasa1999a}). Then, using an It$\hat{o}$
type approximation for the stochastic integral we can obtain
\begin{eqnarray*}
\tilde{l}_n(\theta)&=&\sum^n_{i=1}\frac{\mu(\theta,t_{i-1},X_{t_{i-1}})}{\sigma^2(t_{i-1},X_{t_{i-1}})}(X_{t_i}-X_{t_{i-1}})\\
&-&\frac{1}{2}\sum^n_{i=1}\frac{\mu^2(\theta,t_{i-1},X_{t_{i-1}})}{\sigma^2(t_{i-1},X_{t_{i-1}})}(t_i-t_{i-1}).
\end{eqnarray*}
Thus, the maximizer of $\tilde{l}_n(\theta)$ provides an approximate
maximum likelihood estimate (AMLE). In 1992, Yoshida \cite{Yoshida}
proved that the AMLE is weakly consistent and asymptotically
normally distributed when the diffusion is homogeneous and ergodic.
In ~\cite{Jaya P. N. Bishwal}, Bishwal got the similar result for
the nonhomogeneous case with drift function $\mu(\theta,t, X)=\theta
f(t,X_t)$ for some smooth functions f(t,x). Moreover, he measured
the loss of information using several AMLEs according to different
approximations to $l_T(\theta)$. \\

\subsubsection{Parameter Estimation of Diffusion Terms (and/or Drift Terms) with Discrete
Observation} In previous sections, we always assume one of those
parameters is known and then estimate the other one. In this
section, I want to include the situation when both $\theta$ and
$\vartheta$ are unknown and how to estimate them based on the
discrete observation of the diffusion process at the same time.\\

Suppose we are considering the real valued diffusion process $X_t$
satisfying the following stochastic differential equation
\begin{equation}\label{BMDiffDrift}
dX_t=\mu(\theta, X_t)dt+\sigma(\vartheta, X_t)dB_t.
\end{equation}

Denote the observation times by $t_0=0,t_1, t_2,\ldots, t_{N_T}$,
where $N_T$ is the smallest integer such that $\tau_{N_{T}+1}>T$. In
this section, we mainly consider three cases of estimating
$\beta=(\theta, \vartheta)$, jointly, $\beta=\theta$ with
$\vartheta$ known and $\beta=\vartheta$ with $\theta$ known. In
regular circumstances, the estimate $\hat{\beta}$ converges in
probability to some $\bar{\beta}$ and
$\sqrt{T}(\hat{\beta}-\bar{\beta})$ converges in law to $N(0,
\Omega_{\beta})$ as T tends to infinity, where $\beta_0$ is the true value of the parameter.\\

For simplicity, we set the law of the sampling intervals
$\Delta_n=\tau_n-\tau_{n-1}$ as
\begin{equation}
\Delta=\epsilon \Delta_0,
\end{equation}
where $\Delta_0$ has a given finite distribution and $\epsilon$ is
deterministic.\\
Remark: We are not only studying the case when the sampling interval
is fixed, i.e., $Var[\Delta_0]=0$, but also the continuous
observation case, i.e., $\epsilon=0$ and the random sampling case.\\


Let $h(y_1,y_0,\delta, \beta, \epsilon)$ denote a r-dimensional
vector function which consists of r moment conditions of the
discretized stochastic differential equation \eqref{BMDiffDrift}
(please see \cite{Hansen1982} or \cite{Heyde97} for more details).
Moreover, this function also satisfies
\begin{equation*}
E_{\Delta_n, Y_n, Y_{n-1}}[h(Y_n, Y_{n-1}, \Delta_n, \beta,
\epsilon)]=0,
\end{equation*}
where the expectation is taken with respect to the joint law of
($\Delta_n, Y_n, Y_{n-1}$).\\

By the Law of Large Numbers, $E[h(Y_n, Y_{n-1}, \Delta_n, \beta,
\epsilon)]$ may be estimated by the sample average defined by
\begin{equation}\label{geBMdifanddrisampleaverage}
m_T(\beta)\equiv N^{-1}_T\sum^{N_T-1}_{n=1}h(Y_n, Y_{n-1}, \Delta_n.
\beta, \epsilon).
\end{equation}
Then we can obtain an estimate $\hat{\beta}$ by minimizing a
quadratic function
\begin{equation}\label{geBMdfdrMLF}
Q_T(\beta)\equiv m_T(\beta)'W_T m_T(\beta),
\end{equation}
where $W_T$ is a $r\times r$ positive definite weight matrix and
this method is called Generalized Method of Moments (GMM). In
\cite{Hansen1982}, Hansen proved the strong consistency and
asymptotic normality of GMM estimate, i.e.
\begin{equation*}
\sqrt{N}(\hat{\theta}-\theta)\rightarrow N(0, V),
\end{equation*}
when $\theta=\vartheta$ and $W_T$ satisfied certain conditions.
Mykland used this technique to obtain the closed form for the
asymptotic bias but sacrificed the consistency of the estimate.


\section{Quantifying Uncertainties in SDEs Driven by  Fractional Brownian Motion}

Colored noise, or noise with non-zero correlation in time, are
common in   physical, biological and engineering sciences. One
candidate for modeling colored noise is fractional Brownian motion \cite{DuanLiWang}.

\subsection{Fractional Brownian Motion}

Fractional Brwonian motion (fBM) was   introduced within a Hilbert
space framework by Kolmogorov in 1940 in \cite{Kolmogorov1940},
where it was called \emph{Wiener Helix}. It was further studied by
Yaglom in \cite{Yaglom1958}. The name \emph{fractional Brownian
motion} is due to Mandelbrot and Van Ness, who in 1968 provided in
\cite{Mandelbrot and Van Ness} a stochastic integral representation
of this process in terms of a standard Brownian motion.

\begin{figure}
\begin{minipage}[t]{0.3\linewidth}
\centering
\includegraphics[width=2.2in]{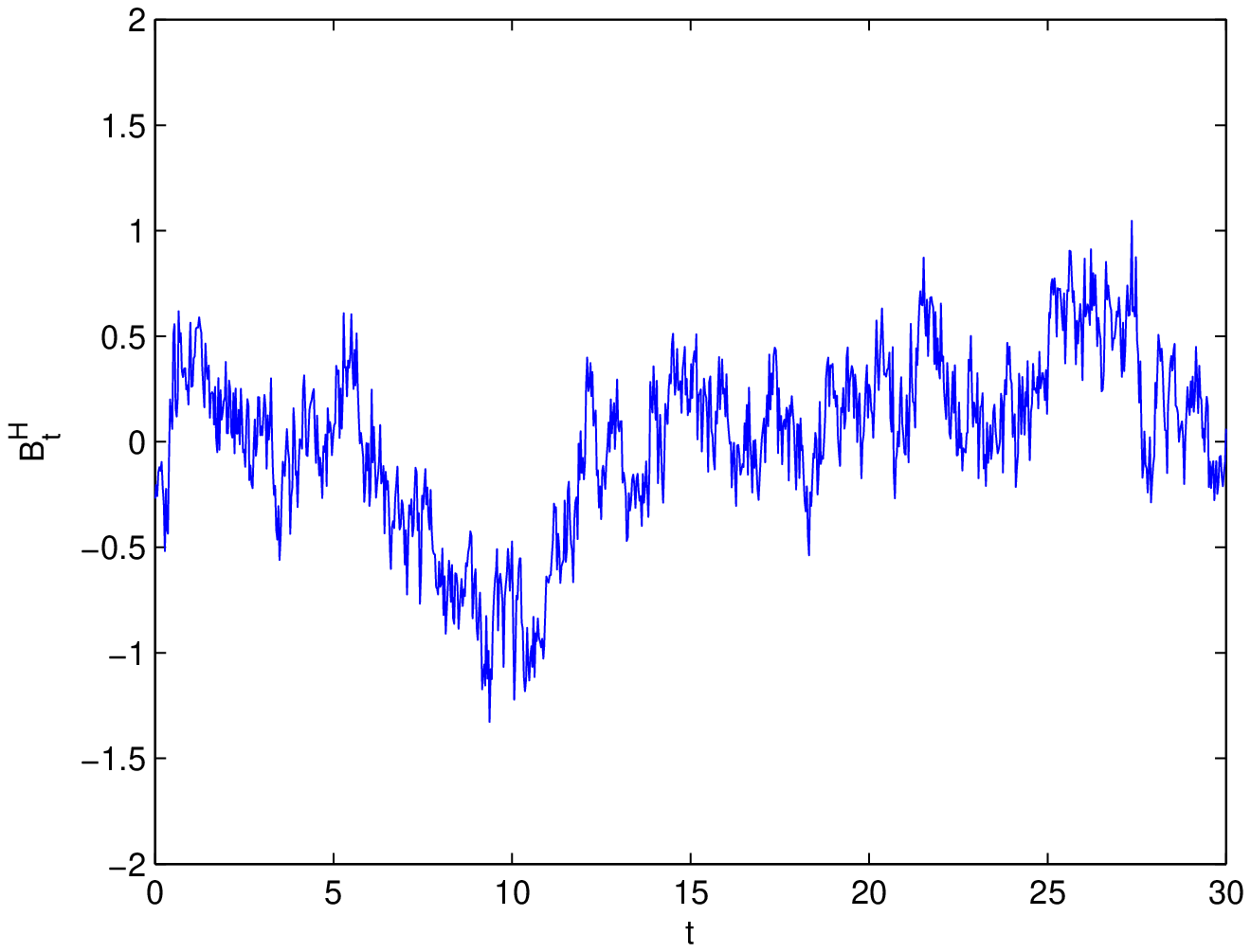}
\end{minipage}%
\begin{minipage}[t]{0.3\linewidth}
\centering
\includegraphics[width=2.2in]{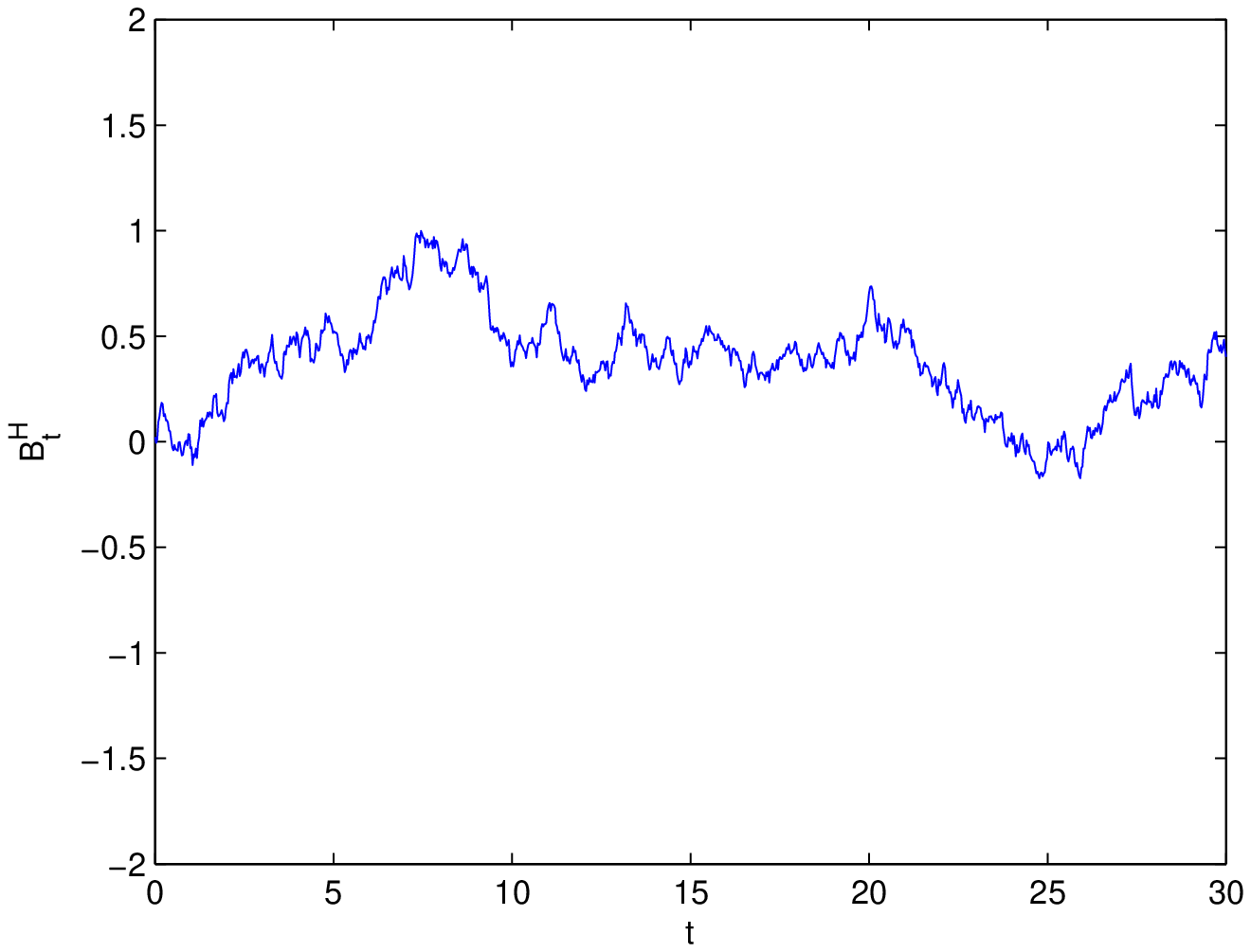}
\end{minipage}
\begin{minipage}[t]{0.3\linewidth}
\centering
\includegraphics[width=2.2in]{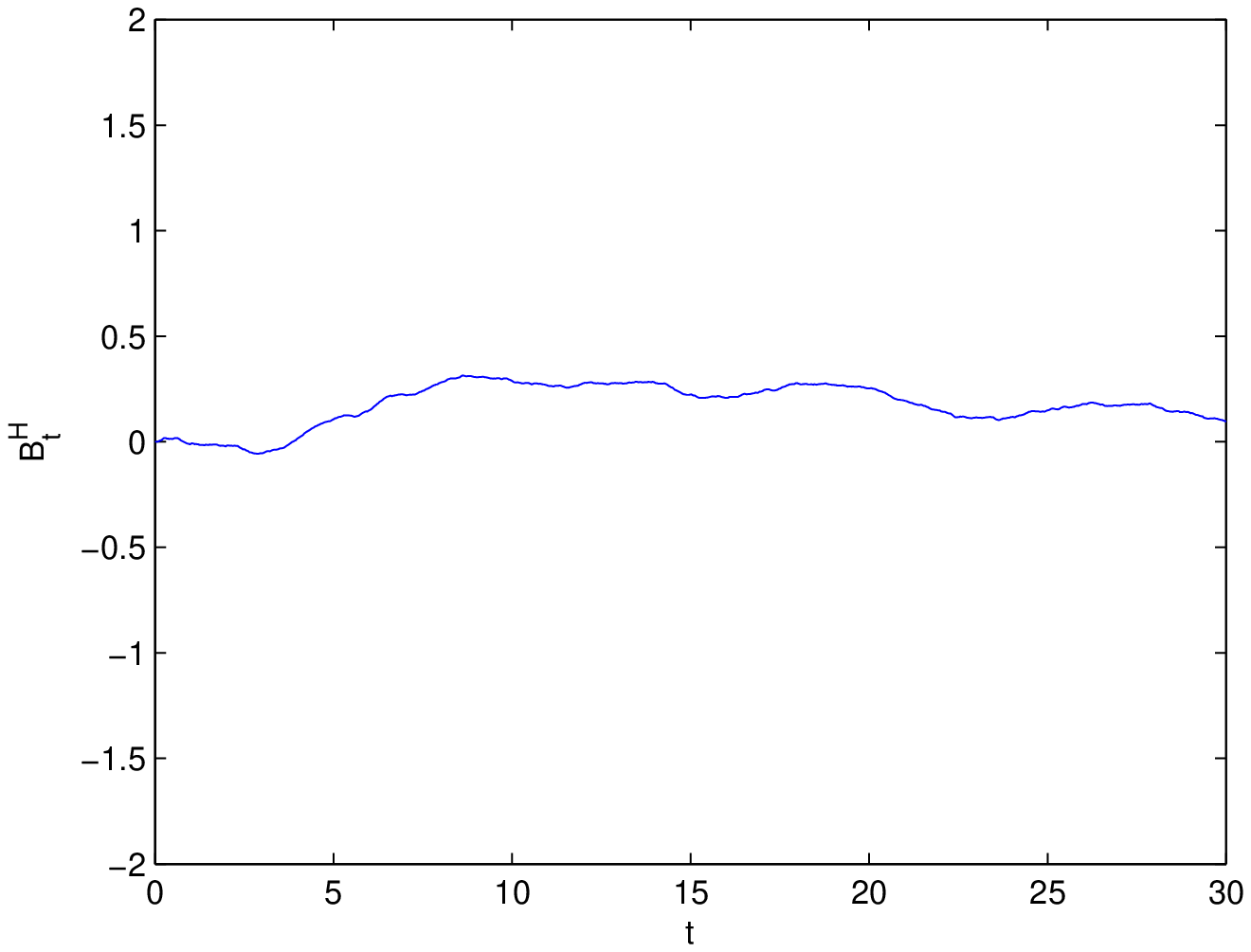}
\end{minipage}
\caption{Three sample paths of fBM with Hurst parameter $H=0.25, 0.5, 0.9$}
\label{fBmPic}
\end{figure}

\begin{definition}[Fractional Brownian motion \cite{Bernt Oksendal}]
Let H be a constant belonging to (0,1). A \emph{fractional Brownian
motion} (fBM) ($B^{H}(t))_{t\geq 0}$ of \emph{Hurst index} H is a
continuous and centered Gaussian process with covariance function
\begin{equation*}
E[B^{H}(t)B^{H}(s)]=\frac{1}{2}(t^{2H}+s^{2H}-|t-s|^{2H}).
\end{equation*}
By the above definition, we see that a standard fBM $B^H$ has the
following properties:
\begin{enumerate}
  \item $B^H(0)=0$ and $E[B^H(t)]=0$ for all $t\geq 0.$
  \item $B^H$ has homogeneous increments, i.e., $B^H(t+s)-B^H(s)$
  has the same law of $B^H(t)$ for $s,t \geq0.$
  \item $B^H$ is a Gaussian process and $E[B^H(t)^2]=t^{2H}, t\leq
  0,$ for all $H\in (0,1)$.
  \item $B^H$ has continuous trajectories.
\end{enumerate}
\end{definition}
Using the method presented in \cite{Coeurjolly1, Coeurjolly2}, we
can simulate sample paths of fractional Brwonian motion with
different Hurst parameters (see Figure \ref{fBmPic}).

For H = 1/2, the fBM is then a standard Brownian motion. Hence, in
this case the increments of the process are independent. On the
contrary, for $H\neq 1/2$ the increments are not independent. More
precisely, by the definition of fBM, we know that the covariance between
$B^H(t+h)-B^H(t)$ and $B^H(s+h)-B^H(s)$ with $s+h\leq t$ and
$t-s=nh$ is
\begin{equation*}
\rho_H(n)=\frac{1}{2}h^{2H}[(n+1)^{2H}+(n-1)^{2H}-2n^{2H}].
\end{equation*}
In particular, we obtain that the two increments of the form
$B^H(t+h)-B^H(t)$ and $B^H(t+2h)-B^H(t+h)$ are positively correlated
for $H>1/2$, while they are negatively correlated for $H<1/2$. In
the first case the process presents an aggregation behavior and this
property can be used in order to describe "cluster" phenomena
(systems with \emph{memory} and \emph{persistence}). In the second
case it can be used to model sequences with \emph{intermittency} and
\emph{antipersistence}.\\

 From the above description, we can get a
general ideal that the Hurst parameter H plays an important role on
how respective fBM behaves. So, it should be considered as an extra
parameter when we estimate others in the coefficients of the SDE
driven by fBM.

Considering the further computation, we would like to introduce one
more useful property of fBM.

\begin{definition}[Self-similarity]
A stochastic process X$=\{X_t, t\in \R\}$ is called
b-\emph{self-similar} or satisfies the property of self-similarity
if for every $a>0$ there exists $b>0$ such that
\begin{equation*}
Law(X_{at},t\geq 0)=Law(a^bX_t, t\geq 0).
\end{equation*}
Note that 'Law=Law' means that the two processes $X_{at}$ and
$a^bX_t$ have the same finite-dimensional distribution functions,
i.e., for every choice $t_0,\ldots, t_n$ in $\R$,
\begin{equation*}
P(X_{at_0}\leq x_0,\ldots, X_{at_n}\leq x_n)=P(a^bX_{t_0}\leq
x_0,\ldots, a^bX_{t_n}\leq x_n).
\end{equation*}
for every $x_0,\ldots, x_n$ in $\R$.
\end{definition}
Since the covariance function of the fBM is homogeneous of order 2H,
we obtain that $B^H$ is a self-similar process with Hurst index H,
i.e., for any constant $a>0$, the processes $B^H(at)$ and
$a^HB^H(t)$ have the same distribution law.

\subsection{How to Estimate Hurst Parameter $H$}
Let's start with the simplest case:
\begin{equation*}
dX_t=dB^H(t),  ~~~~~~~~~~~~~i.e., X_t=B^H(t) ~~~~~~~~~~t\geq 0,
\end{equation*}
where $\{B^H(t), t\geq 0\}$ is a fBM with Hurst parameter $H\in
(0,1)$. Now, our question is how to estimate Hurst parameter H given
observation data $X_0, X_1, \ldots, X_N$. For a parameter estimation
of Hurst parameter H, we need an extra ingredient, fractional
Gaussian noise (fGn).
\begin{definition}[Fractional Gaussian noise] \cite{Gennady
Samorodnitsky2000}\\
Fractional Gaussian noise (fGn) $\{Y_i, i\geq 1\}$ is the increment
of fractional Brownian motion, namely
\begin{equation*}
Y_i=B^H(i+1)-B^H(i),~~~~~~~~~~~~~i\geq 1.
\end{equation*}
\end{definition}
\textbf{Remark:} It is a mean zero, stationary Gaussian time series
whose autocovariance function is given by
\begin{equation*}
\rho(h)=E(Y_iY_{i+h})=\frac{1}{2}\{(h+1)^{2H}-2h^{2H}+|h-1|^{2H}\},~~~h\geq
0.
\end{equation*}
An important point about $\rho(h)$ is
\begin{equation*}
\rho(h)\sim H(2H-1)h^{2H-2},~~~~as~~~h\rightarrow\infty,
\end{equation*}
when $H \neq 1/2$. Since $\rho(h)=0$ for $h\geq 1$ when H=1/2, the
$X_i$'s are white noise in this case. The $X_i$'s, however, are
positively correlated when $\frac{1}{2}<H<1$, and we say that they
display \emph{long-range dependence} (LRD) or \emph{long memory}.

From the expression of fGn, we know it is the same to estimate the
Hurst parameter of fBM as to estimate the Hurst parameter of the
respective fGn. Here, we introduce 4 different methods for measuring
the Hurst parameter. Measurements are given on artificial data and
the results of each method are compared in the end. However, the
measurement techniques used in this paper can only be described
briefly but references to fuller descriptions with mathematical
details are given.


\subsubsection{R/S Method} The R/S method is one of the oldest and
best known techniques for estimating H. It is discussed in detail in
\cite{B. B. Mandelbrot 1969} and \cite{J. Beran}, p.83-87.

For a time series $\{Y_t: t=1,2,\ldots, N \}$ with partial sums
given by $Z(n)=\sum^n_{i=1}Y_i$ and the sample variance given by
\begin{equation*}
S^2(n)=\frac{1}{n-1}\sum^n_{i=1}Y^2_i-\frac{1}{n(n-1)}Z(n)^2,
\end{equation*}
then the R/S statistic, or the \emph{rescaled adjusted range}, is
given by:
\begin{equation*}
\frac{R}{S}(n)=\frac{1}{S(n)}\left[\max_{1\leq t\leq
n}\left(Z(t)-\frac{t}{n}Z(n)\right)-\min_{1\leq t\leq
n}\left(Z(t)-\frac{t}{n}Z(n)\right)\right]
\end{equation*}
For fractional Gaussian noise,
\begin{equation*}
E[R/S(n)]\sim C_Hn^H,
\end{equation*}
as $n\rightarrow \infty$, where $C_H$ is another positive, finite
constant not dependent on n. \\

The procedure to estimate H is therefore as follows. For a time
series of length N, subdivide the series into K blocks with each of
size $n=N/K$. Then, for each lag n, compute $R/S(k_i,n)$, starting
at points $k_i=iN/K+1, i=1,2,\ldots,K-1.$ In this way, a number of
estimates of $R/S(n)$ are obtained for each value of n. For values
of n approaching N, one gets fewer values, as few as 1
when $n\geq N-N/K$.\\

Choosing logarithmically spaced values of n, plot $\log[R/S(k_i,n)]$
versus $\log n$ and get, for each n, several points on the plot.
This plot is sometimes called the pox plot for the R/S statistic.
The parameter H can be estimated by
fitting a line to the points in the pox plot. \\

There are several disadvantages with this technique. Most notably,
there are more estimates of the statistic for low values of n where
the statistic is affected most heavily by short range correlation
behavior. On the other hand, for high values of n there are too few
points for a reliable estimate. The values between these high and
low cut off points should be used to estimate H but, in practice,
often it is the case that widely differing values of H can be found
by this method depending on the high and low cut off points chosen.
To modify the R/S statistic, we can use a weighted sum of
autocovariance instead of the sample variance. Details can be found
in \cite{A. W. Lo.}.

\subsubsection{Aggregated Variance.} Given a time series $\{Y_t:
t=1,2,\ldots,N\}$, divide this into blocks of length m and aggregate
the series over each block.
\begin{equation*}
Y^{(m)}(k):=\frac{1}{m}\sum^{km}_{i=(k-1)m+1}Y_i,~~~~~k=1,2,...,[N/m].
\end{equation*}
We compute its sample variance,
\begin{equation*}
\widehat{VarY^{(m)}}=\frac{1}{N/m}\sum^{N/m}_{k=1}(Y^{(m)}(k)-\overline{Y})^2.
\end{equation*}
where
\begin{equation*}
\overline{Y}=\frac{\sum^N_{i=1}Y_i}{N}.
\end{equation*}
is the sample mean. The sample variance should be asymptotically
proportional to $m^{2H-2}$ for large $N/m$ and m. Then, for
successive values of m, the sample variance of the aggregated series
is plotted versus m on a log-log plot. So we can get the estimate of
H by computing the gradient of that log-log plot. However, jumps in
the mean and slowly decaying trends can severely affect this
statistic. One technique to combat this is to difference the
aggregated variance and work instead with
\begin{equation*}
\widehat{VarY^{(m+1)}}-\widehat{VarY^{(m)}}.
\end{equation*}

\subsubsection{Variance of Residuals.} This method is described in
more detail in \cite{C. K. Peng}. Take the series $\{Y_t:
t=1,2,\ldots,N\}$ and divide it into blocks of length m. Within each
block calculate partial sums:
$Z_k(t)=\sum^{(k-1)m+t}_{i=(k-1)m+1}Y_i$, $k=1\ldots N/m,~t=1\ldots
m$. For each block make a least squares fit to a line $a_k+b_kt$.
Subtract this line from the samples in the block to obtain the
residuals and then calculate their variance
\begin{equation*}
V_k=\frac{1}{m}\sum^{m}_{t=1}(Z_k(t)-a_k-b_kt)^2.
\end{equation*}
The variance of residuals is proportional to $m^{2H}$. For the proof
in the Gaussian case, see \cite{M. S. Taqqu 1997}. This variance of
residuals is computed for each block, and the median (or average) is
computed over the blocks. A log-log plot versus m should follow a
straight line with a slope of 2H.

\subsubsection{Periodogram.} The periodogram is a frequency domain
technique described in \cite{J. Geweke1983}. For a time series
$\{Y_t: t=1,2,\ldots,N\}$, it is defined by
\begin{equation*}
I(\lambda)=\frac{1}{2\pi
N}\left|\sum^N_{j=1}Y_je^{ij\lambda}\right|^2,
\end{equation*}
where $\lambda$ is the frequency. In the finite variance case,
$I(\lambda)$ is an estimate of the spectral density of Y, and a
series with long-range dependence will have a spectral density
proportional to $|\lambda|^{1-2H}$ for frequencies close to the
origin. Therefore, the log-log plot of the periodogram versus the
frequency displays a straight line with a slope of 1-2H.

\subsubsection{Results on Simulated Data} In this subsection, we would like to use
artificial data to check the robustness of above techniques and
compare the result in the end. \\

\begin{figure}
 \includegraphics[height=3in,width=6in]{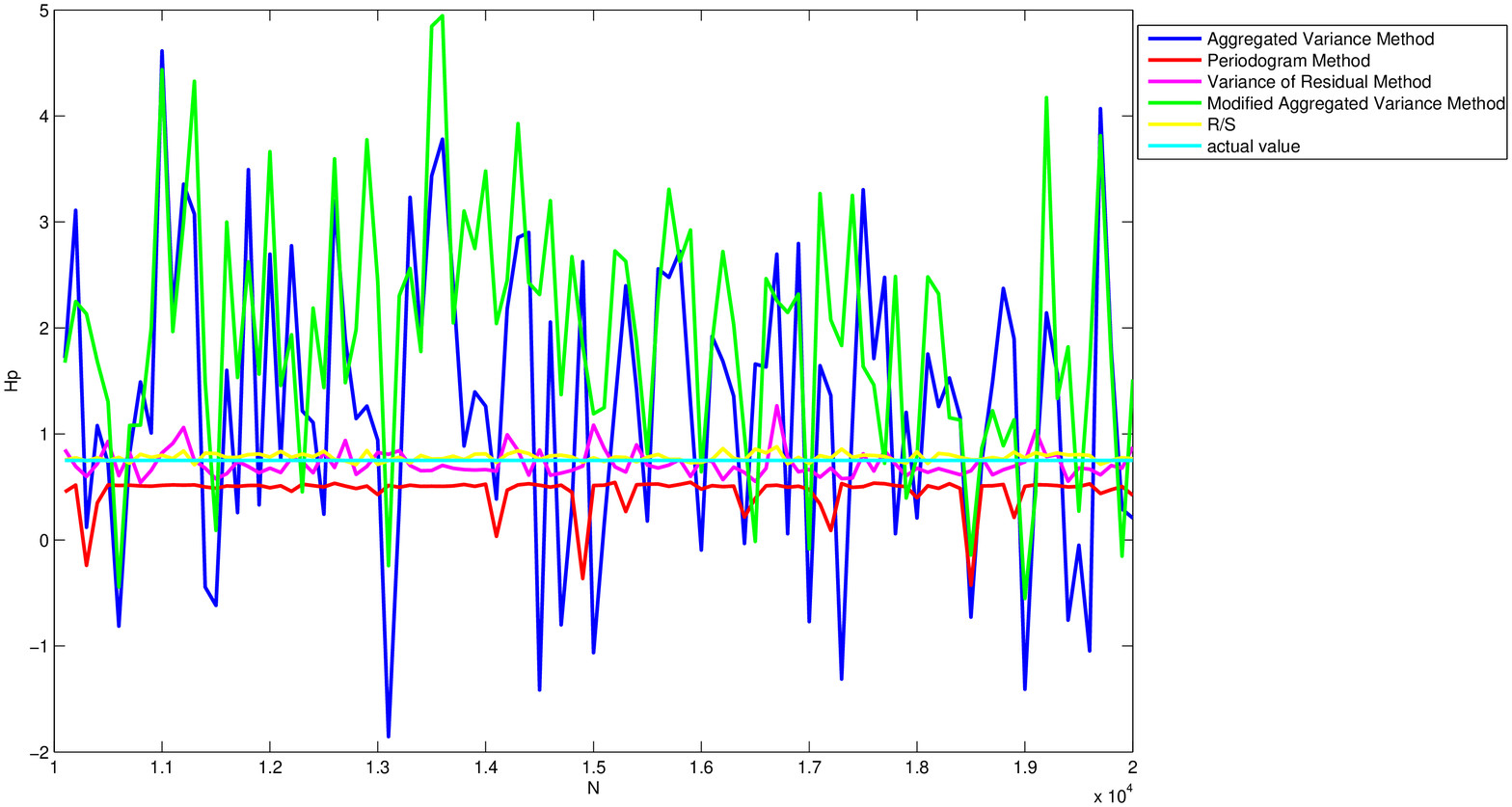}
\caption{Numerical estimation of the Hurst parameter $H$ of fBM:   Actual
value $H=0.65$} \label{Hurst1}
\end{figure}

\begin{figure}
 \includegraphics[height=3in,width=6in]{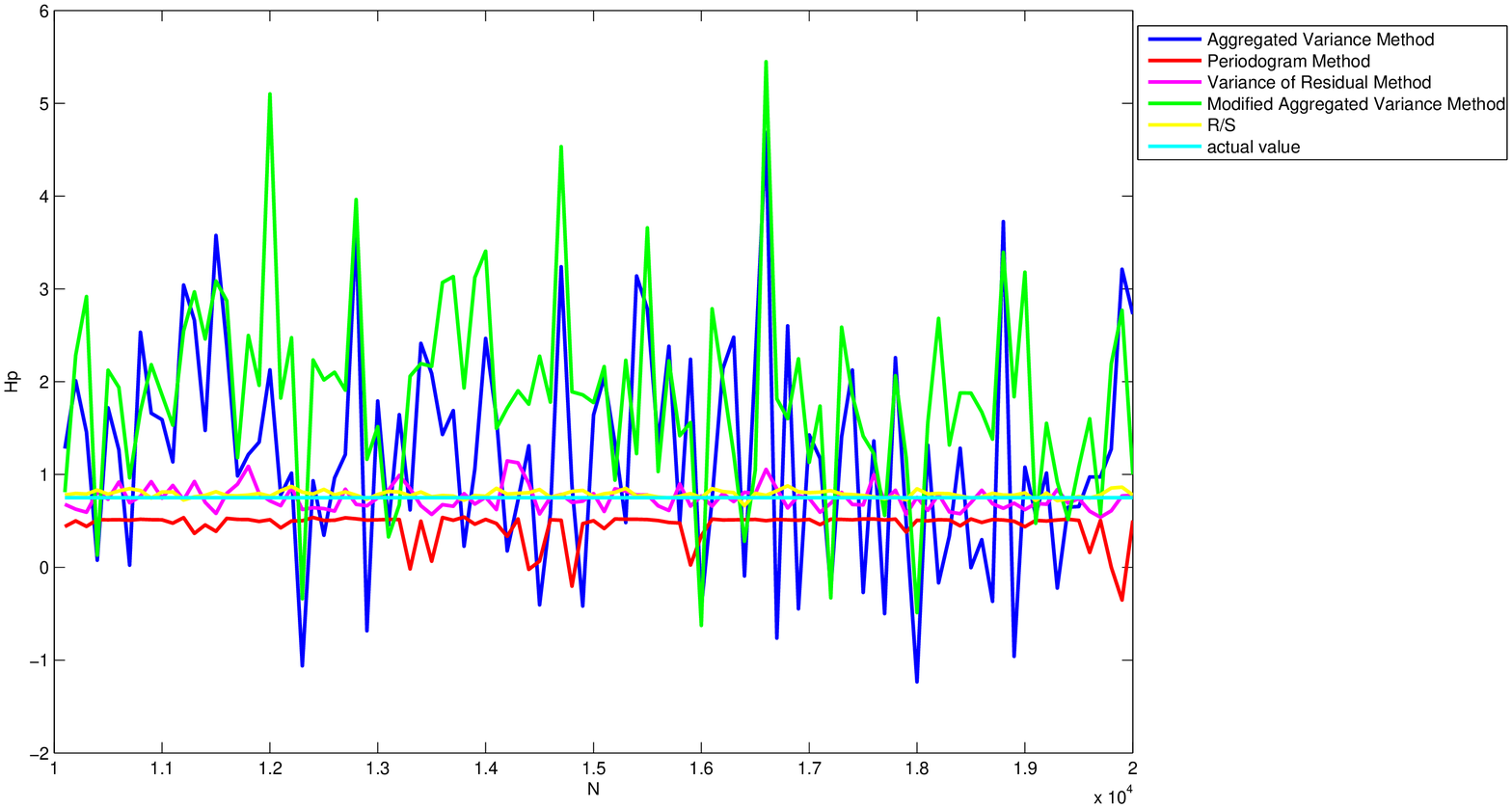}\\
\caption{Numerical estimation of the Hurst parameter $H$ of fBM:   Actual
value $H=0.95$} \label{Hurst2}
\end{figure}

For each of the simulation methods chosen, traces have been
generated. Each trace is 10,000 points of data. Hurst parameters of
0.65 and 0.95 have been chosen to represent a low and a high level
of long-range dependence in data. From the Figure \ref{Hurst1} and
Figure \ref{Hurst2}, we can see that the Variance of Residual Method
and R/S have the most accurate result. The Modified Aggregated
Variance Method improved a little bit over the original one, but
both of them still fluctuate too much.\\

\subsection{How to Estimate Parameters in SDEs Driven by fBM}
After we discuss how to estimate the Hurst parameter of a series of
artificial fBM data, now we want to concern how to estimate the
parameters of the linear/nonlinear stochastic differential
equation(s) driven by fBM. The coefficients in the stochastic
differential equation could be deterministic or random, linear or
nonlinear. No general results are available. So some
specific statistical results will be discussed below according to what
kind of specified models we deal with.\\

\subsubsection{Preparation} The main difficulty in dealing with a fBm is
that it is not a semimartingale when $H\neq \frac{1}{2}$ and hence
the results from the classical stochastic integration theory for
semimartingales can not be applied. So, we would like to introduce
the following integral transformation which can transform fBM to
martingale firstly and it will be a key point in our development
below. For $0<s<t\leq T$, denote
\begin{eqnarray}\label{martingale}
k_H(t,s)&=&\kappa^{-1}_Hs^{(1/2)-H}(t-s)^{(1/2)-H},\\
\kappa_H&=&2H\Gamma(3/2 -H)\Gamma(H+1/2),\\
w^H_t&=&\lambda^{-1}_Ht^{2-2H};~~~~~\lambda_H=\frac{2H\Gamma(3-2H)\Gamma(H+1/2)}{\Gamma(3/2-H)},  \\
M^H_t&=&\int^t_0 k_H(t,s)dB^H_s.
\end{eqnarray}
Then the process $M^H$ is a Gaussian martingale (see \cite{Le
Breton} and \cite{Norros1999}), called the \emph{fundamental
martingale} with variance function $w^H$.

\subsubsection{Parameter Estimation for a Fractional Langevin Equation}

Suppose $\{X_t, ~t\geq 0\}$ satisfies the following stochastic
differential equation
\begin{equation*}\label{fBMLangevin}
X_t=\theta\int^t_0X_s ds+\s B^H_t;~~~~~~~~~0\leq t\leq T,
\end{equation*}
where $\t$ and $\s^2$ are unknown constant parameters, $B^H_t$ is a
fBM with Hurst parameter $H\in [1/2, 1]$. \\

Denote the process Z=$(Z_t, t\in [0, T])$ by
\begin{equation}\label{Z}
Z_t=\int^t_0 k_H(t,s)dX_s.
\end{equation}
Then we can prove that process Z is a semimartingale associated to X
with following decomposition (see \cite{Kleptsyna2000})
  \begin{equation}\label{semimartingal expression}
  Z_t=\theta \int^t_0 Q(s) dw^H_s+\sigma M^H_t,
  \end{equation}
  where
  \begin{equation}\label{Q}
  Q(t)=\frac{d}{dw^H_t}\int^t_0 k_H(t,s)X(s)d s,
  \end{equation}
  and $M^H_t$is the Gaussian martingale defined by
  (17).
From the representation \eqref{semimartingal expression}, we know
the quadratic variation of Z on the interval [0, t] is nothing but
\begin{equation*}
\langle Z\rangle_t=\sigma^2w^H_t,~~~~~~a.s.
\end{equation*}
Hence the parameter $\sigma^2$ can be obtained by
\begin{equation*}
[w^H_t]^{-1}\lim_n\sum_i
\left[Z_{t^{n}_{i+1}}-Z_{t^{n}_{i}}\right]^2=\sigma^2,~~~~~~a.s.
\end{equation*}
where $t^{n}_i$ is an appropriate partition of [0,t] such that
$\sup_i|t^{n}_{i+1}-t^{n}_{i}|\rightarrow 0$ as $n \rightarrow
\infty$. So, the variance parameter can be computed with
probability 1 on any finite time interval.\\

As for  the parameter $\theta$, by applying the Girsanov type formula
for fBM which is proved in \cite{Kleptsyna2000}, we can define the
following maximum likelihood estimate of $\theta$ based on the
observation on the interval [0, t] by
\begin{equation}\label{MLE-FOU}
\theta_T=\left\{\int^T_0 Q^2(s)dw^H_s\right\}^{-1}\int^T_0 Q(s)dZ_s,
\end{equation}
where processes Q, Z and $w^H_t$ are defined by \eqref{Q}, \eqref{Z}
and (16), respectively. For this estimate, strong consistency is
proven and explicit formulas for the asymptotic bias and mean square
error are derived by Kleptsyna and Le Breton \cite{Kleptsyna2002}.

Remarks: \begin{itemize}
  \item When $H=1/2$, since $Q=Z=X$ and $d\omega_s^{1/2}=ds$, the formula
\eqref{MLE-FOU} reduces to the result of   \cite{Lipster1977} for
an usual Ornstein-Uhlenbeck process.
  \item For an arbitrary $H\in [1/2,1]$, we could derive the following
alternative expression for $\theta_T$:
\begin{equation*}
\theta_T=\left\{2\int^T_0
Q^2(s)dw^H_s\right\}^{-1}\left\{\frac{\lambda_H}{2-2H}Z_T\int^T_0s^{2H-1}dZ_s-t\right\}.
\end{equation*}
\end{itemize}

\begin{example}
Consider a special   Ornstein-Uhlenbeck model
\begin{equation*}
dX_t=\theta X_tdt+  2 dB^H_t.
\end{equation*}
Then, according to the above approximation scheme, we can  numerically estimate
  $\theta=1$   and the results are
shown in Figure \ref{fBM332}.

\begin{figure}
 \includegraphics[height=2.9in,width=6.4in]{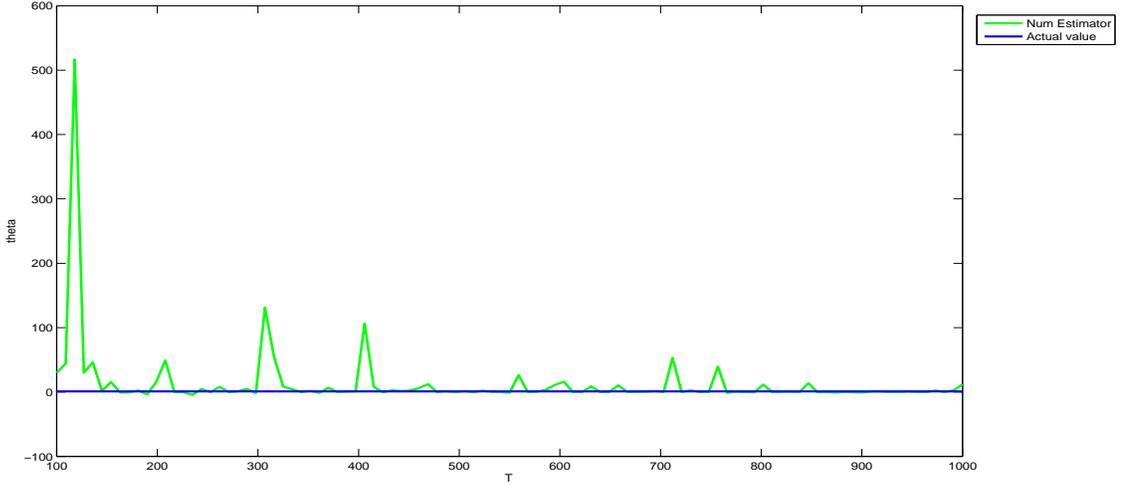}
\caption{Numerical estimation of drift parameter $\theta$ in $
dX_t=\theta X_tdt+ 2 dB^H_t $  with Hurst parameter $H=0.75$:
Actual value $\theta_0=1$} \label{fBM332}
\end{figure}
\end{example}

\subsubsection{Parameter Estimation in Linear Deterministic Regression}
Suppose $X_t$ satisfies the following stochastic differential
equation
\begin{equation*}
X_t=\theta \int^t_0 A(s) ds+ \int^t_0 C(s) dB^H_s,~~~~~~~~~0\leq
t\leq T,
\end{equation*}
where A and C are deterministic measurable functions on [0,T],
$B^H_t$ is a fBM with Hurst parameter $H\in [1/2, 1]$. \\

Let $q_t$ be defined by
\begin{equation*}
q_t=\frac{d}{dw^H_t}\int^t_0 k_H(t,s)\frac{A}{C}(s)d s,
\end{equation*}
where $w^H_t$ and $k_H(t,s)$ are defined by (16) and
\eqref{martingale}. Then, from Theorem 3 in \cite{Kleptsyna2000}, we
obtain the maximum likelihood estimate of $\theta$ defined by
\begin{equation*}\label{MLE-LDR}
\theta_T=\left\{\int^T_0 q^2_t dw^H_t\right\}^{-1} \int^T_0 q_t
dZ_t,
\end{equation*}
where $Z_t$ is defined by \eqref{Z}. \\

Remark: This result can be extended to an arbitrary H in (0,1) (see
\cite{Le Breton}) and $\theta_T$ is also the best linear unbiased
estimate of $\theta$.

\begin{example}
Consider a special Linear Deterministic Regression
\begin{equation*}
dX_t=-\theta dt+t dB^H_t.
\end{equation*}
Then, using the above estimate, we can do numerical simulation with
result shown in Figure \ref{fBM333}.
\begin{figure}
 \includegraphics[height=3in,width=6in]{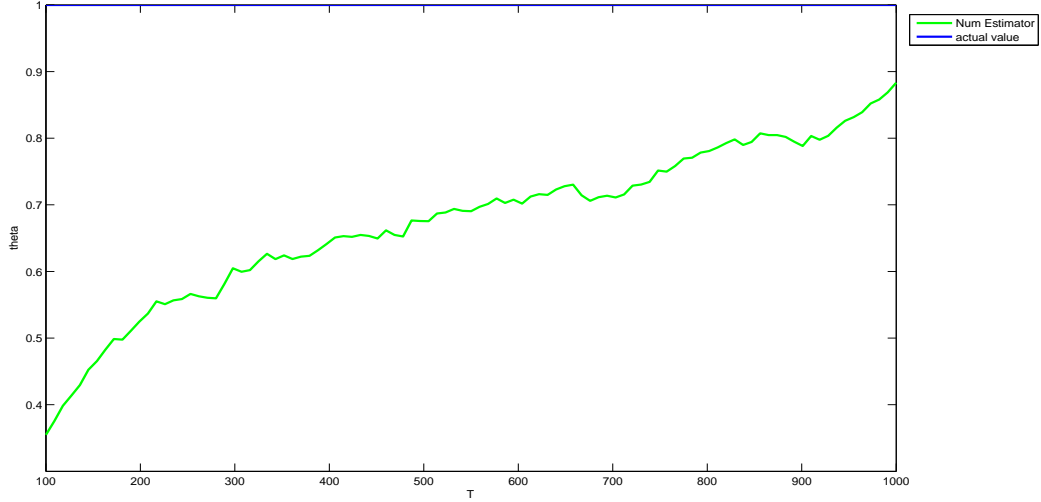}
\caption{Numerical estimation of drift parameter $\theta$ in a Linear
Deterministic Regression  $dX_t= -\theta  dt+t dB^H_t$ with  Hurst parameter $H=0.75$:   Actual value
  $\theta$=1}
\label{fBM333}
\end{figure}
\end{example}

\subsubsection{Parameter Estimation in Linear Random Regression}
Let us consider a stochastic differential equation
\begin{equation*}
\d X(t)=[A(t,X(t))+\theta C(t,X(t))]dt+\sigma(t)dB^H_t,~~~t\geq 0,
\end{equation*}
where $B=\{B^H_t, t\geq 0\}$ is a fractional Brownian motion with
Hurst parameter H and $\sigma(t)$ is a positive nonvanishing
function  on $[0, \infty)$. According to \cite{B. L. S. Prakasa
Rao2003}, the Maximum Likelihood Estimate $\hat{\t}_T$ of $\t$ is
given by
\begin{equation*}\label{MLE-LRR}
\t_T=\frac{\int^T_0 J_2(t)dZ_t+\int^T_0
J_1(t)J_2(t)dw^H_t}{\int^T_0J^2_2(t)dw^H_t},
\end{equation*}
where the processes $Z_t, J_1, J_2$ are defined by
\begin{eqnarray*}
Z_t&=&\int^t_0\frac{k_H(t,s)}{\sigma(s)}dX_s, t\geq 0,\\
J_1(t)&=&\frac{d}{dw^H_t}\int^t_0
k_H(t,s)\frac{A(s,X(s))}{\sigma(s)}ds,~~~~~J_2(t)=\frac{d}{dw^H_t}\int^t_0
k_H(t,s)\frac{C(s,X(s))}{\sigma(s)}ds, \\
\end{eqnarray*}
and $w^H_t$, $k_H(t,s)$ are defined by (16) and (14). Also in the
same paper, they proved that $\t_T$ is strongly consistent for the
true value $\t$.

\begin{example}
Consider a special Linear Random   Regression
\begin{equation*}
dX_t=(t+\theta X_t) dt+t dB^H_t.
\end{equation*}
A numerical estimation of the parameter $\theta$ is
 shown in Figure \ref{fBM334}.
\begin{figure}
 \includegraphics[height=3in,width=6in]{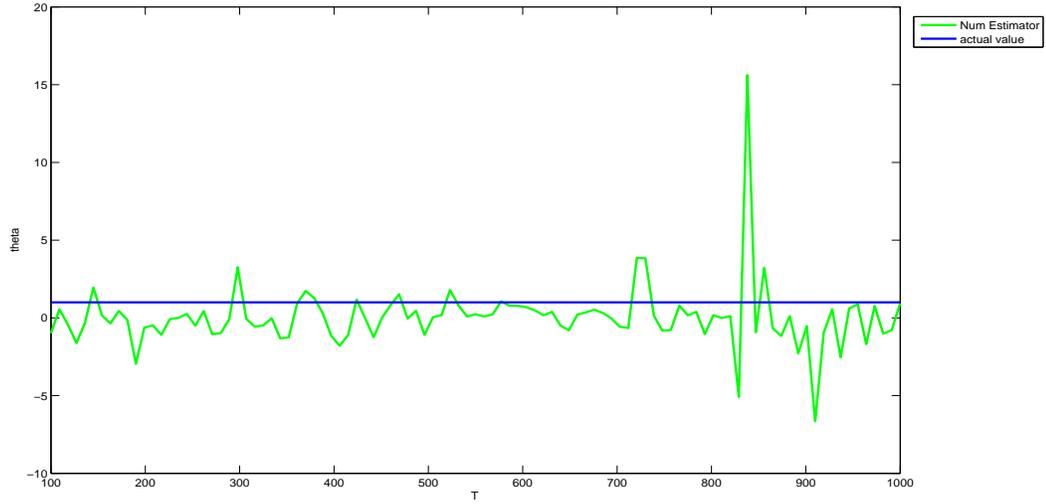}
\caption{Numerical estimation of drift parameter $\theta$ in a Linear
Random Regression  $dX_t=(t+\theta X_t) dt+t dB^H_t$ with  Hurst parameter $H=0.75$:   Actual value
  $\theta$=1}
\label{fBM334}
\end{figure}
\end{example}

\section{Parameter Estimation for SDE Driven by $\alpha$-Stable L\'evy Motion}

  Brownian motion, as a Gaussian process, has been widely used to
model fluctuations in   engineering and science. For a particle in
Brownian motion, its sample paths are continuous in time almost
surely (i.e., no jumps), its mean square displacement increases
linearly in time (i.e., normal diffusion), and  its probability
density function decays exponentially in space (i.e., light tail or
exponential relaxation) \cite{Oksendal}.  However some complex phenomena
involve non-Gaussian fluctuations, with   properties such as
anomalous diffusion (mean square displacement is a nonlinear power
law of time) \cite{BD90} and heavy tail (non-exponential relaxation)
\cite{Yon96}. For instance, it has been argued that diffusion in a
case of geophysical turbulence \cite{Shlesinger} is anomalous.
Loosely speaking, the diffusion process consists of a series of
``pauses", when the particle is trapped by a coherent structure, and
``flights" or ``jumps" or other extreme events, when the particle
moves in a jet flow. Moreover, anomalous electrical transport
properties have been observed in some amorphous materials such as
insulators, semiconductors and polymers, where transient current is
asymptotically a power law function of time \cite{SSB91, Herrchen}.
Finally, some paleoclimatic data \cite{Dit} indicates heavy tail
distributions and some DNA data \cite{Shlesinger} shows long range
power law decay for spatial correlation. L\'evy motions are thought
to be appropriate models for non-Gaussian processes with jumps
\cite{Sato-99}. Here we consider a special non-Gaussian process, the
$\alpha$-stable L\'evy motion, which arise in many complex systems \cite{Woy}.

\subsection{$\alpha$-Stable L\'evy Motion}
There are several reasons for using a stable distribution to model a
 fluctuation process in a dynamical system. Firstly, there are   theoretical reasons for expecting a
non-Gaussian stable model, e.g. hitting times for a Brownian motion
yielding a L\'evy distribution, and reflection off a rotating mirror
yielding a Cauchy distribution. The second reason is the Generalized
Central Limit Theorem which states that the only possible
non-trivial limit of normalized sums of i.i.d. terms is stable. The
third argument for modeling with stable distributions is empirical:
Many large data sets exhibit heavy tails and skewness. In this section,
we consider one-dimensional $\alpha$-stable distributions defined as
follows.
\begin{definition}(\cite{Aleksander}, Chapter 2.4) The \textbf{Characteristic Function} $\phi(u)$ of an \textbf{$\alpha$-stable random
variable} is given by
\begin{equation}
\phi(u)=\exp((-\sigma^{\alpha})|u|^{\a}\{1-i\b sgn(u) \tan(\a
\pi/2)\}+i\m u)
\end{equation}
where $\a \in (0,1)\cup(1,2),~\b\in[-1,1],\s\in \R_+,~\m\in\R$, and
by
\begin{equation}
\phi(u)=\exp(-\sigma|u|\{1+i\b\frac{2}{\pi} sgn(u)\log(|u|)\}+i\m u)
\end{equation}
when $\alpha=1$, it gives a very well-known symmetric Cauchy
distribution and
\begin{equation}
\phi(u)=\exp(-\frac{1}{2}\s|u|^2+i\m u),
\end{equation}
when $\alpha=2$, it gives the well-known Gaussian
distribution.
\end{definition}

For the random variable X distributed according to the rule
described above we use the notation $X\thicksim S_{\a}(\s,\b,\m)$.
Especially when $\m=\b=0$, i.e., X is a symmetric $\a$-stable random
variable, we will denote it as $X\thicksim S\a S$.\\

Also, from above definition, it is easy to see that the full stable
class is characterized by four parameters, usually designated
$\alpha, \b, \s,$ and $\mu$. The shift parameter $\mu$ simply shifts
the distribution to the left or right. The scale parameter $\s$
compresses or extends the distribution about $\mu$ in proportion to
$\s$ which means, if the variable x has the stable distribution
$X\thicksim S_{\a}(\s,\b,\m)$, the transformed variable
$z=(x-\mu)/\s$ will have the same shaped distribution, but with
location parameter 0 and scale parameter 1. The two remaining
parameters completely determine the distribution's shape. The
characteristic exponent $\alpha$ lies in the range $(0,2]$ and
determines the rate at which the tails of the distribution taper
off. When $\alpha=2$, a normal distribution results. When
$\alpha<2$, the variance is infinite. When $\alpha>1$, the mean of
the distribution exists and is equal to $\mu$. However, when
$\alpha\leq 1$, the tails are so heavy that even the mean does not
exist. The fourth parameter $\b$ determines the skewness of the
distribution and lies in the range [-1,1].

Now let us introduce   $\a$-stable L$\e$vy motions.
\begin{definition}($\a$-stable L$\e$vy motion   \cite{Aleksander})\\
A stochastic process $\{X(t): t\geq 0\}$ is called the (standard)
$\a$-stable L$\e$vy motion if
\begin{enumerate}
  \item X(0)=0 a.s.;
  \item $\{X(t): t\geq 0\}$ has independent increments;
  \item X(t)-X(s)$\thicksim S_{\a}((t-s)^{1 / \a},\b,0)$ for any $0\leq s<t<\infty$.
\end{enumerate}
\end{definition}
So, from the third condition,  we can simulate all
$\a$-stable L$\e$vy motion if we know how to simulate $X\thicksim
S_{\a}(\s,\b,0)$.  Especially, it is enough to simulate $X\thicksim
S_{\a}(\s,0,0)$ if we want to get the trajectories of symmetric
$\a$-stable L$\e$vy motions.\\

We recall  an important property of $\a$-Stable random
variables giving us the following   result: It is enough to know
how to simulate $X\thicksim S_{\a}(1,0,0)$ in order to get any
$X\thicksim S_{\a}(\s,0,0), \forall \s\in \R^+$.
\begin{prop}If we have $X_1,X_2\thicksim
S_{\a}(\s,\b,\m)$ and A, B are real positive constants and C is a
real constant, then
\begin{equation*}
AX_1+BX_2+C\thicksim S_{\a}(\s(A^{\a}+B^{\a})^{1/\a},\b,
\m(A^{\a}+B^{\a})^{1/\a}+C)\\
\end{equation*}
\end{prop}

\begin{prop}Let $X\sim S_{\alpha}(\s,\b,0)$, with $0<\alpha<2$, Then
$E|X|^p<\infty$ for any $0<p<\alpha$, $E|X|^p=\infty$ for any $p\geq
\alpha.$
\end{prop}

\begin{figure}
\begin{minipage}[t]{0.3\linewidth}
\centering
\includegraphics[width=2.2in]{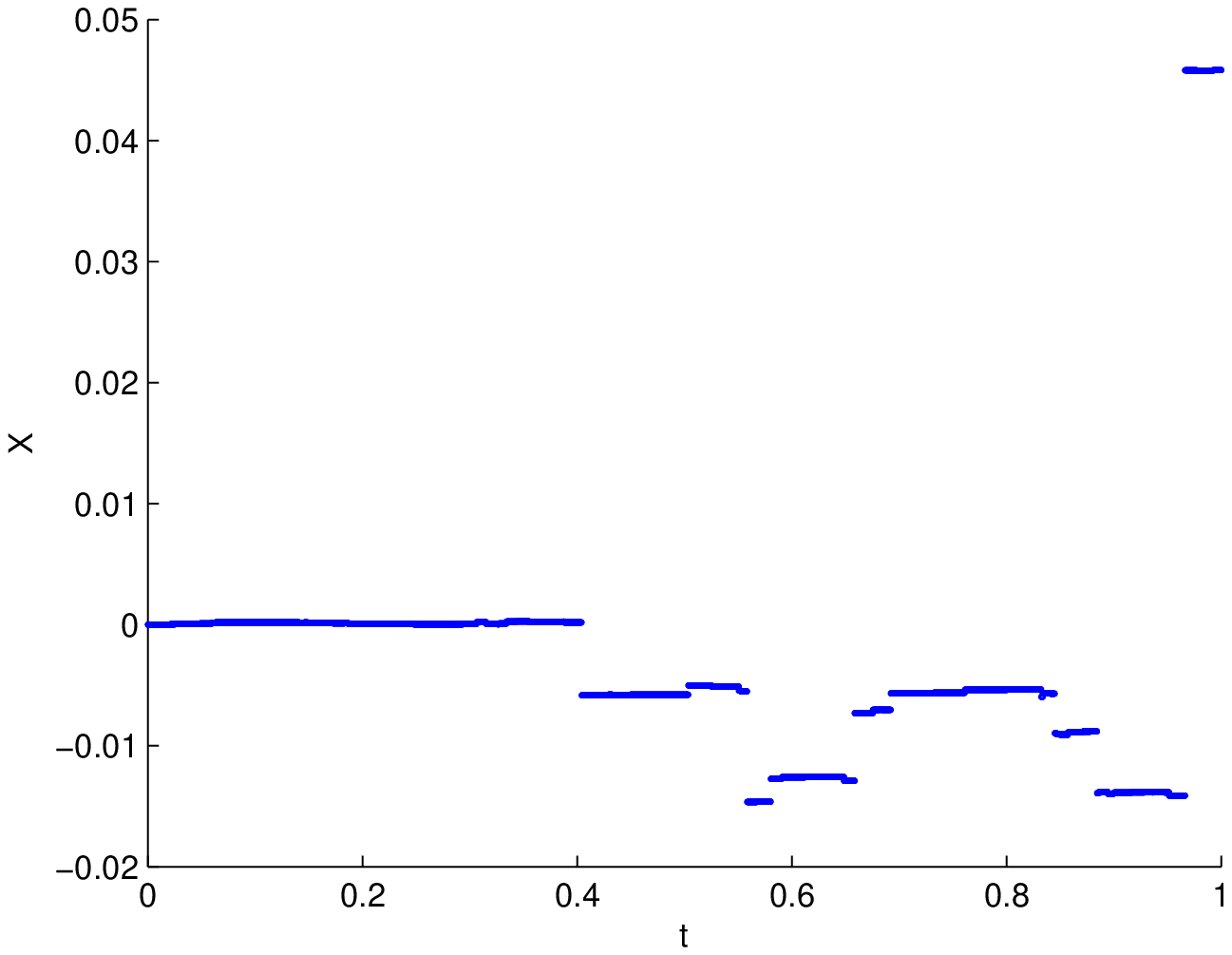}
\end{minipage}%
\begin{minipage}[t]{0.3\linewidth}
\centering
\includegraphics[width=2.2in]{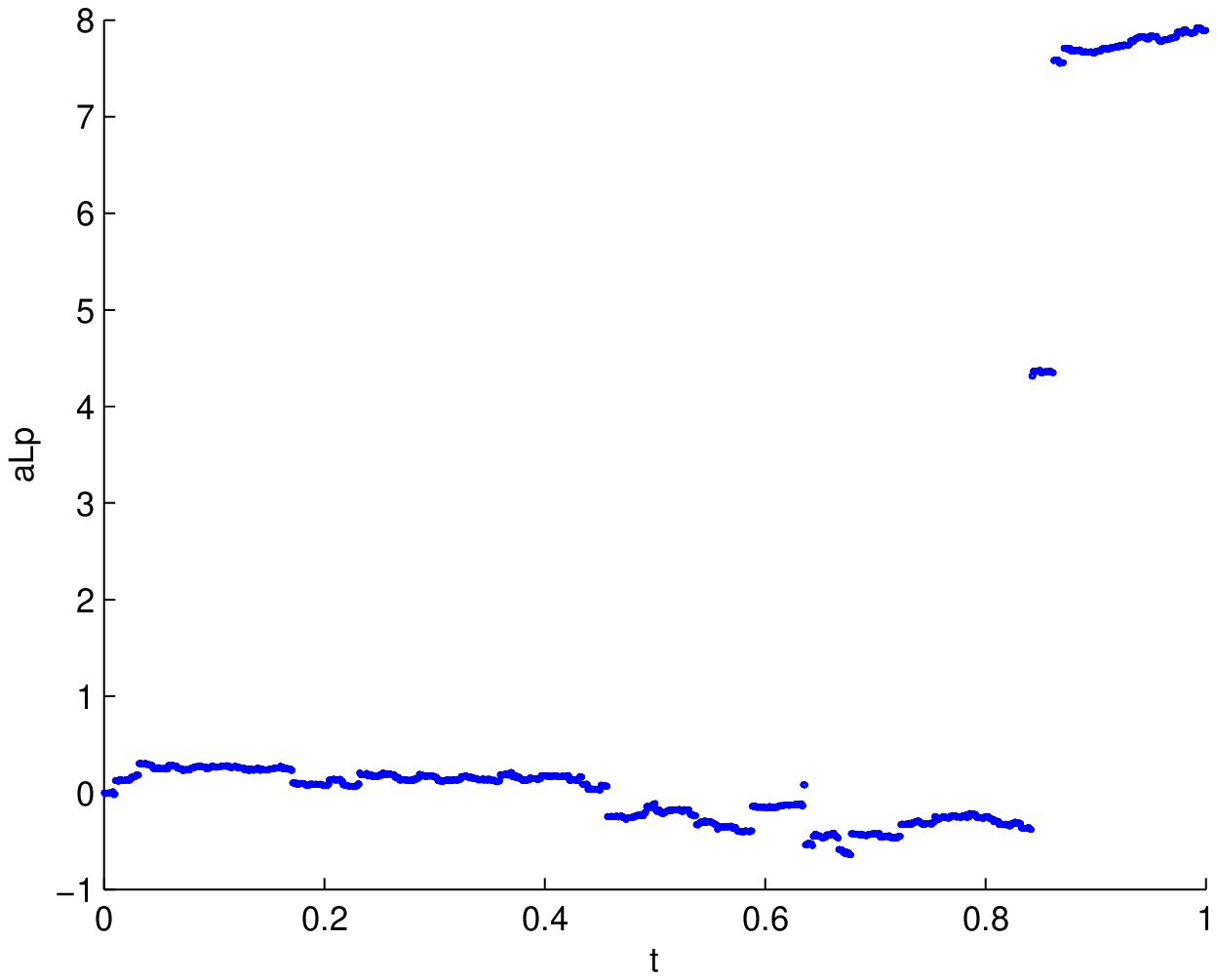}
\end{minipage}
\begin{minipage}[t]{0.3\linewidth}
\centering
\includegraphics[width=2.2in]{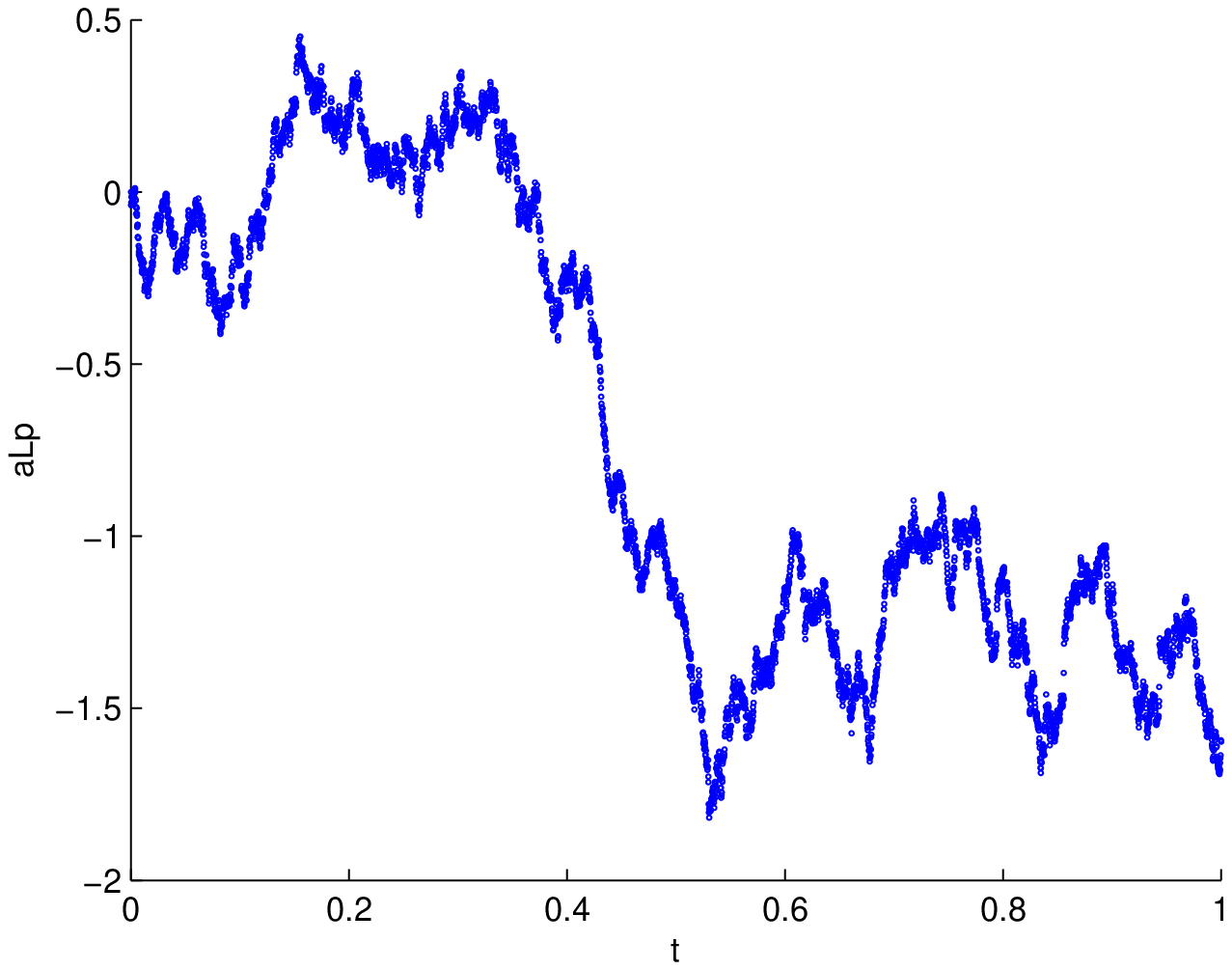}
\end{minipage}
\caption{Three sample paths of symmetric $\alpha-$stable L\'evy
motion with $\alpha=0.4, 1.2, 1.9$, respectively} \label{Levypic}
\end{figure}

Figure \ref{Levypic} shows sample paths of the $\alpha$-stable
L\'evy motion with different $\alpha$.
\\

As we can see in Figure \ref{Levypic}, the bigger the parameter
$\alpha$ is, the more the path looks like Brownian motion. Generally
speaking, when we deal with concrete data, we have to choose
$\alpha$-stable processes very carefully to get the best estimation.
We now discuss how to estimate $\alpha$.

\subsection{How to Estimate the Characteristic Exponent $\alpha$}
Five different methods about how to estimate the
characteristic exponent $\alpha$ of $\alpha-$stable distribution are
considered: Characteristic Function Method(CFM), Quantile Method,
Maximum Likelihood Method, Extreme Value Method and Moment Method.
As in the last section, measurements are given on artificial data
and the results of each method are compared in the end of this
section.

\subsubsection{Characteristic Function Method}

Since  $\alpha-$stable distributions are uniquely determined by
their Characteristic Function (CF), it is natural to consider how to
estimate parameter by studying their CF.   Press \cite{Press 1972}
introduced a parameter estimation method based on CF, which gets
estimations of parameters by minimizing differences between values
of sample CF and the real ones. But this method is
only applicable to standard distributions.\\

Another method which uses the linearity of logarithm of CF was
developed by Koutrouvelis \cite{Koutrouvelis} and it can be applied
to general $\alpha$-stable cases. This method is denoted as Kou-CFM.
The
  idea is as follows: On the one hand, taking the
logarithm of real part of CF gives
\begin{equation*}
\ln[-Re(\phi(u
))]=\alpha \ln|u|+\alpha \ln \sigma.
\end{equation*}
On the other hand, the sample characteristic function of $\phi(\t)$
is given by $\Phi(\t)=(\sum^N_{k=1}e^{i\t y_k})$ where $y_k$'s are $n$
independent observations. In \cite{Koutrouvelis}, a regression
technique is applied to gain estimates for all parameters of a
observed $\alpha$ stable distribution. In \cite{Kogon}, Kogon
improved this method by replacing a linear regression fit by a
linear least square fit which gave a more accurate estimate and its
computational complexity became lower.

\subsubsection{Quantile Method}
Quantiles are points taken at regular intervals from the cumulative
distribution function of a random variable. Suppose we have $n$
independent symmetric $\alpha$-stable random variables with the
stable distribution $S_{\alpha}(\s, \b, \mu)$, whose parameters are
to be estimated. Let $x_p$ be the p-th quantile, so that
$S_{\alpha}(x_p; \s, \b, \mu)=p$. Let $\hat{x_p}$ be the
corresponding sample quantile, then $\hat{x_p}$ is a consistent
estimate of $x_p$. \\

In 1971, Fama and Roll \cite{Fama 1971} discovered that, for some large p (for
example, p=0.95),
\begin{equation*}
\hat{z}_p=\frac{\hat{x}_p-\hat{x}_{1-p}}{2\sigma}=(0.827)\frac{\hat{x}_p-\hat{x}_{1-p}}{\hat{x}_{0.72}-\hat{x}_{0.28}}
\end{equation*}
is an estimate of the p-quantile of the standardized symmetric
stable distribution with exponent $\alpha$. According to this, they
proposed a estimate (QM) for symmetric $\alpha$-stable
distributions. However, the serious disadvantage of this method is
that
its estimations are asymptotically biased.\\
Later on,   McCulloch \cite{McCulloch} improved and extended
this result to general $\alpha$-stable distributions, denoted as
McCulloch-QM. Firstly, he defined
\begin{eqnarray*}
v_{\alpha}&=&(x_{-0.95}-x_{-0.05})/(x_{-0.75}-x_{-0.25})\\
v_{\beta}&=&(x_{-0.95}+x_{-0.05}-2x_{0.5})/(x_{-0.95}-x_{-0.05})
\end{eqnarray*}
and let $\hat{v}_{\alpha}$ and $\hat{v}_{\beta}$ be the
corresponding sample value:
\begin{eqnarray*}
\hat{v}_{\alpha}&=&(\hat{x}_{-0.95}-\hat{x}_{-0.05})/(\hat{x}_{-0.75}-\hat{x}_{-0.25})\\
\hat{v}_{\beta}&=&(\hat{x}_{-0.95}+\hat{x}_{-0.05}-2\hat{x}_{0.5})/(\hat{x}_{-0.95}-\hat{x}_{-0.05})
\end{eqnarray*}
which are the consistent estimates of the index $v_{\alpha}$ and
$v_{\beta}$. Then, he illustrated that estimates of $\alpha$ can be
expressed by a function of $\hat{v}_{\alpha}$ and $\hat{v}_{\beta}$.
Compared with QM, McCulloch-QM could get consistent and unbiased
estimations for the general $\alpha$-stable distribution, and extend
the estimation range of parameter $\alpha$ to $0.6\leq \alpha \leq
2$. Despite its computational simplicity, this method has a number
of drawbacks, such as, there are no analytical expressions for the
value of the fraction, and the evaluation of the tables implies that
it is highly dependent on the value of $\alpha$ in a nonlinear way.
This technique does not provide any closed-form solutions.


\subsubsection{Extreme Value Method}
In 1996, based on asymptotic extreme value theory, order statistics
and fractional lower order moments,   Tsihrintzis and
  Nikias   \cite{Tsihrintzis1996} proposed a new
estimate which can be computed fast for symmetric $\alpha$ stable
distribution from a set of i.i.d. observations. Five years later,
Kuruoglu \cite{Kuruoglu} extended it to the general $\alpha$ stable
distributions. The general idea of this method is as follows. Given a
data series $\{X_i: i=1,2, \ldots, N\}$, divide this into L
nonoverlapping blocks of length K such that $K=N/L$. Then the
logarithms of the maximum and minimum samples of each segment are
computed as follows
\begin{eqnarray*}
\overline{Y_l}&=&\log(\max\{X_{lK-K+i}| i=1,2,\ldots, K\}),\\
\underline{Y_l}&=&\log(-\min\{X_{lK-K+i}| i=1,2,\ldots, K\}).
\end{eqnarray*}
The sample means and variances of $\overline{Y_l}$ and
$\underline{Y_l}$ are calculated as
\begin{eqnarray*}
\overline{Y}&=&\frac{1}{L}\sum^L_{l=1}\overline{Y_l},~~~~~\overline{s^2}=\frac{1}{L-1}\sum^L_{l=1}(\overline{Y_l}-\overline{Y})^2, \\
\underline{Y}&=&\frac{1}{L}\sum^L_{l=1}\underline{Y_l},~~~~~\underline{s^2}=\frac{1}{L-1}\sum^L_{l=1}(\underline{Y_l}-\underline{Y})^2.
\end{eqnarray*}
Finally, an estimate for $\alpha$ is given by sample variance as
follows
\begin{equation*}
\hat{\alpha}=\frac{\pi}{2\sqrt{6}}\left(\frac{1}{\overline{s}}+
\frac{1}{\underline{s}}\right).
\end{equation*}
Even though the accuracy and computational complexity decrease, there is
now a closed form for the block size which means we have to look-up
table to determine the segment size K.

\subsubsection{Moment Estimation Method}
Another way to estimate parameters of the general $\alpha$-stable
distribution is the Logarithmic Moments Methods which was also
introduced by Kuruoglu \cite{Kuruoglu}. The advantage of this method
relative to the Fractional Lower Order Method is that it does not
require the inversion of a sinc function or the choice of a moment
exponent p. The main feature is that the estimate of $\alpha$ can be
expressed by a function of the second-order moment of the skewed
process, i.e.
\begin{equation*}
\hat{\alpha}=\left(\frac{L_2}{\psi_1}-\frac{1}{2}\right)^{-1/2},
\end{equation*}
where $\psi_1=\frac{\pi^2}{6}$ and, for any $X\sim
S_{\alpha}(\s,\b,0)$, $L_2$ is defined as follows
\begin{equation*}
L_2=E[(\log
|X|-E[\log|X|])^2]=\psi_1\left(\frac{1}{2}+\frac{1}{\alpha^2}\right)-\frac{\theta^2}{\alpha^2}.
\end{equation*}

\subsubsection{Results on Simulated Data}  In this subsection, we
would like to use artificial data to check the robustness of the
above techniques and compare the results. \\

For each of the simulation methods chosen,  estimates of
$\alpha$ have been generated respectively and each trace is 1,000
points of data. Characteristic exponents of 0.95 and 1.70 have been
chosen to represent a low and a high level of the rate at which the
tails of the distribution taper off.

\begin{figure}
 \includegraphics[height=3in,width=6in]{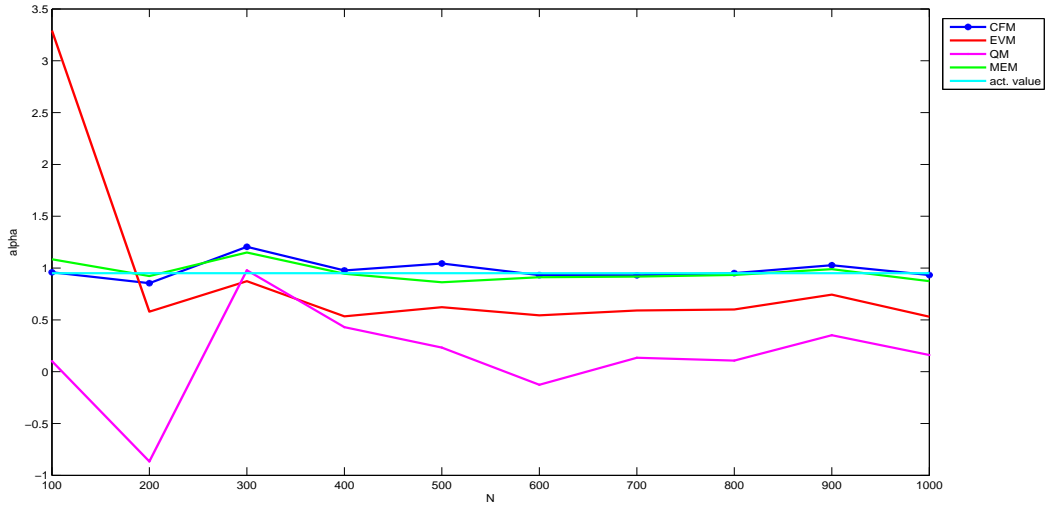}
\caption{Numerical estimation of the characteristic exponent
$\alpha$ in the $\alpha-$stable L\'evy motion $L_t^\a$ by 4
different methods:  Actual value $\a=0.95$} \label{Levy421}
\end{figure}

\begin{figure}
 \includegraphics[height=3in,width=6in]{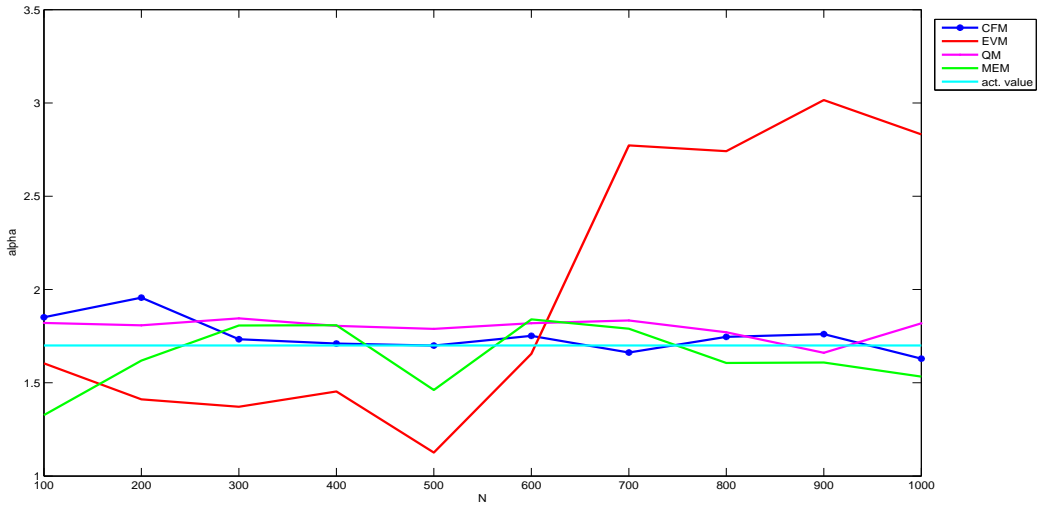}\\
\caption{Numerical estimation of the characteristic exponent
$\alpha$ in the $\alpha-$stable L\'evy motion $L_t^\a$ by 4
different methods:  Actual value $\a=1.70$} \label{Levy422}
\end{figure}

From the Figures \ref{Levy421} and \ref{Levy422}, we can see that
the Characteristic Function Method and the Moment Estimate Method
have the most accurate result. The Quantile Method behaved a little
better than Extreme Value Method, but both of them still fluctuate
too much when $\alpha$ is small. As to the convergence, we can see
that all the methods get closer and closer to the real value when
the points of data increase except for the Extreme Value Method.

\subsection{How to Estimate Parameters in SDEs Driven by L\'evy Motion}
After we discussed how to estimate the characteristic exponent of
$\alpha$-stable L\'evy motions, now we consider how to estimate the
parameters in  stochastic differential equations driven by general
L\'evy motion. Just as what we discussed about fBM, no general
results about the parameter estimation for general cases are
available at this time. Some special results will be listed below
for different equations.

We consider parameter estimation of the L\'evy-driven stationary
Ornstein-Uhlenbeck process.  Recently,  Brockdwell, Davis and  Yang
\cite{Brockwell2007} studied parameter estimation problems for
  L\'evy-driven Langevin equation (whose solution is called an Ornstein-Uhlenceck process)  based on
observations made at uniformly and closely-spaced times. The idea is
to obtain a highly efficient estimate of the
L\'evy-driven Ornstein-Uhlenceck process coefficient by estimating
the corresponding coefficient of
the sampled process. The main feature is discussed below.\\

Consider a stochastic differential equation   driven by the L\'evy motion
$\{L(t), t\geq 0\}$
\begin{equation*}
dY(t)=- \theta Y(t)dt+ \sigma dL(t).
\end{equation*}
When $ L(t)$ is Brownian motion, the solution of
above equation can be expressed as
\begin{equation}\label{solution}
Y(t)=e^{- \theta t}Y(0)+\sigma\int^t_0e^{-\theta(t-u)}dL(t).
\end{equation}
For any second-order driving L\'evy motion, the process $\{Y(t)\}$
can be defined in the same way, and if $\{L(t)\}$ is non-decreasing,
$\{Y(t)\}$ can also be defined pathwise as a Riemann-Stieltjes
integral by \eqref{solution}. For the convenience of the simulation,
we rewrite solution as follows
\begin{equation}\label{solutioniteration}
Y(t)=e^{-\theta (t-s)}Y(s)+\sigma\int^t_s
e^{-\theta (t-u)}dL(u),~~for~all~t>s\geq 0.
\end{equation}
Now we  collect all information corresponding to the sampled
process in order to get the estimate. Set $t=nh$ and $s=(n-1)h$ in
equation \eqref{solutioniteration}. Then the sampled process
$\{Y^{(h)}_n, n=0,1,2,\ldots\}$ (or the discrete-time AR(1) process)
satisfies
\begin{equation*}
Y^{(h)}_n=\phi Y^{(h)}_{n-1}+Z_n,
\end{equation*}
where
\begin{equation*}
\phi=e^{-\theta h},~~~and~~~Z_n=\sigma\int^{nh}_{(n-1)h}e^{-\theta
(nh-u)}dL(u).
\end{equation*}
Then, using the highly efficient Davis-McCormick estimate of $\phi$,
namely
\begin{equation*}
\hat{\phi}^{(h)}_N=min_{1\leq n\leq
N}\frac{Y^{(h)}_n}{Y^{(h)}_{n-1}},
\end{equation*}
we can get the estimate of $\theta$ and $\sigma$ as follows
\begin{eqnarray*}
\hat{\theta}^{(h)}_N&=&-h^{-1}\log \hat{\phi}^{(h)}_N,  \\
\hat{\sigma}^{(2)}_N&=&\frac{2
\hat{\theta}^{(h)}_N}{N}\sum^N_{i=0}(Y^{(h)}_i-\overline{Y}^{(h)}_N)^2.
\end{eqnarray*}
\begin{example}
Consider a    L\'evy-driven   Ornstein-Uhlenbeck
process   satisfying the following SDE
\begin{equation}
dX_t=-   X_t dt+  \sigma  dL^\a_t.
\end{equation}
A numerical estimation of the diffusion parameter $\sigma$
is shown in Figure \ref{LevyOU1}.

\begin{figure}
 \includegraphics[height=3in,width=6in]{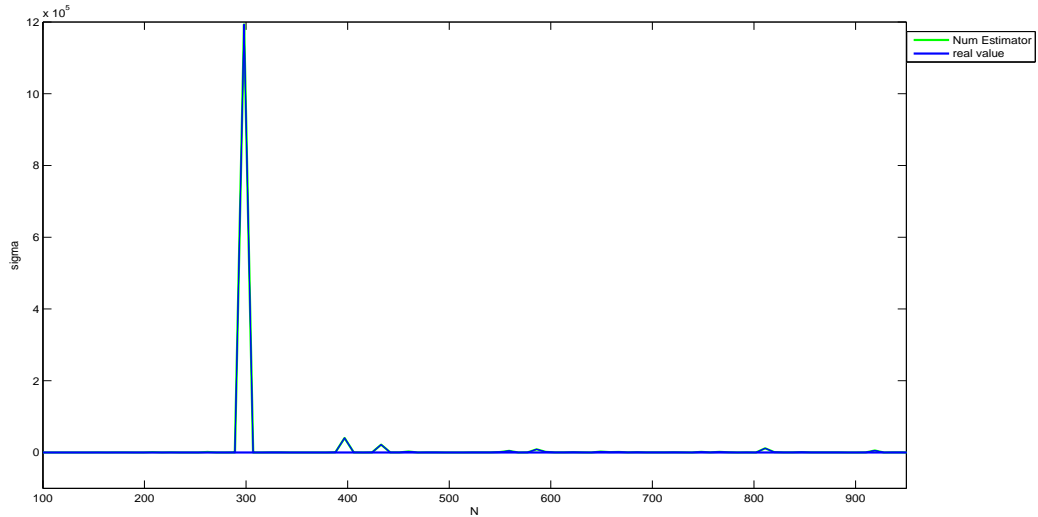}
\caption{Numerical estimation of the diffusion parameter $\sigma$ in
L\'evy-driven Ornstein-Uhlenbeck process  $dX_t=-  X_t dt+  \sigma
dL^\a_t$   with     $\alpha=0.95$: Actual value $\sigma=2$}
\label{LevyOU1}
\end{figure}

\end{example}


\begin{thebibliography}{909}



\bibitem{Ait2002}
Y. A$\ddot{i}$t-Sahalia (2002),
\newblock{ Maximum-likelihood
estimation of discretely-sampled diffusions: a closed-form
approximation approach}, \emph{Econometrica} \textbf{70}, 223-262.

\bibitem{Mykland2004}
Y. A$\ddot{i}$t-Sahalia and  P. A. Mykland (2004), \newblock{
Estimators of diffusions with randomly spaced discrete observations:
a general theory}, \emph{The Annals of Statistics} \textbf{32}(5),
2186-2222


 \bibitem{Ait-Sahalia2003}Y. A$\ddot{i}$t-Sahalia and P. A. Mykland
 (2003),
 \newblock {The effects of random and discrete sampling when estimating continuous-time
 diffusions},
 \newblock \emph{Econometrica} \textbf{71}(2), 483-549.

\bibitem{wu}
S. Albeverrio, B. R\"{u}diger  and J. L. Wu (2000), \newblock{
Invariant measures and symmetry property of L\'{e}vy type
operators}, \emph{Potential Analysis} \textbf{13}, 147-168.

\bibitem{Alizadeh2002} S. Alizadeh, M. W. Brandt and F. X. Diebold
(2002),
 \newblock {Range-based estimation of stochastic volatility models},
 \emph{The Journal of Finance} \textbf{57}(3), 1047-1091.


\bibitem{Applebaum}
D. Applebaum (2009),
\newblock {L\'evy Processes and Stochastic Calculus}, 2nd edition, Cambridge University Press, UK.






\bibitem{Arnold} L. Arnold (1998),
 \newblock {Random Dynamical Systems},
 \newblock Springer, New York.



\bibitem{Barn} O. E. Barndorff-Nielsen, T. Mikosch and S. I. Resnick (Eds.) (2001),
L\'evy Processes: Theory and Applications, Birkh\"auser, Boston.



\bibitem{Bender} C. Bender (2003),  \newblock{An Ito formula for generalized
functionals of a Fractional Brownian motion with arbitrary Hurst
parameter}, \emph{Stoch. Proc. Appl.} \textbf{104}, 81-106.


\bibitem{J. Beran}J. Beran (1994),
 \newblock {Statistics for Long-Memory Processes},
 \newblock Chapman and Hall.

\bibitem{Bertoin-98}
J.~Bertoin (1998), \newblock{L\'evy Processes}, Cambridge University
Press, Cambridge, U.K..

\bibitem{Billingsley1961} P. Billingsley (1961), \newblock{Statistical Inference for Markov Processes},
 Chicago University Press, Chicago.

\bibitem{Jaya P. N. Bishwal} J. P. N. Bishwal (2007),
\newblock{Parameter Estimation in Stochastic Differential Equations},
Springer, New York.

\bibitem{BlomkerStani} D.  Blomker and S. Maier-Paape (2003),
\newblock{Pattern formation below criticality forced by noise}, \emph{Z.
Angew. Math. Phys.} \textbf{54}(1), 1--25.



\bibitem{BD90} J. P. Bouchaud and A. Georges (1990), \newblock{Anomalous diffusion
in disordered media: Statistic mechanics, models and physical
applications}, \emph{Phys. Repts} \textbf{195}, 127-293.


\bibitem{Brockwell2007}P. J. Brockwell, R. A. Davis,
and Y. Yang (2007),
 \newblock {Estimation for nonnegative L\'evy-driven Ornstein-Uhlenbeck
 processes},
 \newblock \emph{J. Appl. Probab.} \textbf{44}(4), 977-989.


 \bibitem{CarLanRob01}
T.~Caraballo, J.~Langa  and J.~C. Robinson (2001), A stochastic
pitchfork bifurcation in a reaction-diffusion
 equation, \emph{R. Soc. Lond. Proc. Ser. A} \textbf{457}, 2041--2061.

\bibitem{ChenBH}B. Chen (2009), Stochastic dynamics of water vapor in the climate
system, Ph.D. Thesis, Illinois Institute of Technology, Chicago,
USA.

\bibitem{ChenDuan} B. Chen and J. Duan (2009),
Stochastic quantification of missing mechanisms in dynamical
systems, In ``Recent Development in Stochastic Dynamics and
Stochastic Analysis", Interdisciplinary Math, Sci. \textbf{8},
67-76.

\bibitem{Chronopoulou2009} A. Chronopoulou and F. Viens (2009),
Hurst index estimation for self-similar processes with long-memory.
In ``Recent Development in Stochastic Dynamics and Stochastic
Analysis", J. Duan, S. Luo and C. Wang (Eds.), World Scientific.

\bibitem{Nualart2006}J. M. Corcuera, D. Nualart, and J. H. C. Woerner
(2006), Power variation of some integral fractional processes,
Bernoulli, \textbf{12}, 713-735.


\bibitem{Nualart2007}J. M. Corcuera, D. Nualart and  J. H. C. Woerner (2007),
A functional central limit theorem for the realized power variation
of integrated stable process, \emph{Stochastic Analysis and
Applications} \textbf{25}, 169-186.


\bibitem{Coeurjolly1} J. Coeurjolly (2001),
 \newblock {Estimating the parameters of the Fractional Brwonian motion by discrete variations of its sample
 paths},
 \emph{Statistical Inference for Stochastic Processes} \textbf{4}, 199-227.

\bibitem{Coeurjolly2}J. Coeurjolly (2000):
 \newblock {Simulation and identification of the Fractional Brwonian motion: a bibliographical and comparative study},
 \emph{Journal of Statistical Software, American Statistical Association} \textbf{5}(07).

\bibitem{Crauel} H. Crauel  and F. Flandoli (1998),
Additive noise destroys a pitchfork bifurcation, \emph{Journal of
Dynamics and Differential Equations} \textbf{10}, 259-274.

\bibitem{Dacumha1986} D. Dacunha-Castelle adn D. Florens-Zmirou (1986),
Estimation of the coefficients of a diffusion from discrete
observations, \textbf{19}, 263-284.

\bibitem{DaPrato} G. Da Prato and J. Zabczyk (1992), {Stochastic Equations
in Infinite Dimensions}, Cambridge University Press.


 \bibitem{Davis2000}M. Davis (2001),
 \newblock {Pricing weather derivatives by marginal value}, \emph{Quantitative Finance} \textbf{1}(3),
305-308.

\bibitem{Dit}P. D. Ditlevsen (1999),
 \newblock {Observation of $\alpha-$stable
noise induced millennial climate changes from an ice record},
\emph{Geophys. Res. Lett.} \textbf{26}, 1441-1444.

\bibitem{Doob1953}J. L. Doob (1953),
 \newblock {\em Stochastic Processes,}
John Wiley, New York.

\bibitem{Du2009}A. Du and J. Duan (2009), \newblock{A stochastic approach for parameterizing
  unresolved scales in a system with memory}, \emph{Journal of Algorithms \&
Computational Technology} \textbf{3}, 393-405.


\bibitem{Duan2009}J. Duan (2009),
\newblock{Stochastic modeling of unresolved scales in complex
systems}, \emph{Frontiers of Math.} in China, \textbf{4}, 425-436.

\bibitem{Duan2009}J. Duan (2009),
\newblock{Predictability in spatially extended systems
  with model uncertainty I \& II}, \emph{Engineering  Simulation} \textbf{2}, 17-32 \& \textbf{3} 21-35.

\bibitem{Duan2009}J. Duan (2009),
\newblock{Predictability in nonlinear dynamical systems
 with model uncertainty}, \emph{Stochastic Physics and Climate
 Modeling}, T. N. Palmer and P. Williams (eds.), Cambridge Univ.
 Press, pp.105-132.



 \bibitem{DuanKanSchm} J. Duan, X. Kan and B. Schmalfuss (2009),
 Canonical sample spaces for stochastic dynamical systems,\emph{
 In ``Perspectives in Mathematical Sciences", Interdisciplinary Math. Sci.} \textbf{9}, 53-70.

\bibitem{DuanLiWang} J. Duan, C. Li and and X. Wang (2009),
Modeling colored noise by Fractional Brownian motion,
\emph{Interdisciplinary Math. Sci.} \textbf{8}, 119-130.

\bibitem{DuanBalu}J. Duan  and B. Nadiga (2007),
Stochastic parameterization of large Eddy simulation of geophysical
flows, \emph{Proc. American Math. Soc.} \textbf{135}, 1187-1196.

\bibitem{Dohnal1987}G. Dohnal (1987),
On estimating the diffusion coefficient, \emph{J. Appl. Prob.}
\textbf{24}, 105-114.

\bibitem{Elerian2000}O. Elerian, S. Chib and N. Shephard (2001),
 \newblock {Likelihood inference for discretely observed non-linear
 diffusions},
\emph{Econometrica} \textbf{69}(4), 959-993.


\bibitem{Eraker1998} B. Eraker (2001),
 \newblock {MCMC analysis of diffusion models with application to finance},
 \emph{Journal of Business and Economic Statistics} \textbf{19}(2),
 177-191.

\bibitem{Fama 1971} E. F. Fama, R. Roll (1971),
\newblock{Parameter estimates for symmetric stable distribution}
\emph{Journal of the American Statistical Association}, \textbf{66},
331-338.

\bibitem{Florens1989} D. Florens-Zmirou (1989),
\newblock{Approximate discrete-time schemes for statistics of diffusion
processes},
\newblock {Statistics} \textbf{20}, 547-557.

\bibitem{Gardiner} C. W. Gardiner (1985), Handbook of Stochastic
Methods, Second Ed., Springer, New York.

\bibitem{Gar} J. Garcia-Ojalvo and J. M. Sancho (1999),
{Noise in Spatially Extended Systems}, Springer-Verlag, 1999.

\bibitem{Valentine and Jacod 1993}V. Genon-Catalot and J. Jacod
(1993),
 \newblock {On the estimation of the diffusion coefficient for
multi-dimensional diffusion processes}, \emph{Annales de l'I. H.
P.}, section B, tome 29, 1993.

\bibitem{Genon} V. Genon-Catalot and J. Jacod (1993),
\newblock{On the estimation of the diffusion coefficient for
multidimensional diffusion processes}, \emph{Ann. Inst. Henri
Poincar\'e, Probabilit\'es et Statistiques.} \textbf{29}, 119-151.

\bibitem{Genon1994} V. Genon-Catalot and J. Jacod (1994),
\newblock{On the estimation of the diffusion coefficient for
diffusion processes}, \emph{J. Statist.} \textbf{21}, 193-221.

 \bibitem{Genon-Catalot1999} V. Genon-Catalot, T. Jeantheau and C. Laredo
 (1999),
 \newblock {Parameter estimation for discretely observed stochastic volatility
 models},
 \newblock \emph{Bernoulli} 5(5), 855-872.

\bibitem{J. Geweke1983} J. Geweke and S. Porter-Hudak (1983),
 \newblock {The estimation and application of long memory time series
 models},
 \newblock \emph{Time Ser. Anal.} \textbf{4}, 221-238.


\bibitem{Hanggi3} P. Hanggi and P. Jung (1995), Colored noise in dynamical
systems, \emph{Advances in Chem. Phys.} \textbf{89}, 239-326.

\bibitem{Hansen1982} L. P. Hansen (1982), Large sample properties of generalized method of moments
estimators,
 Econometrica \textbf{63}, 767-804.

\bibitem{Imkeller09} C.\ Hein, P.\ Imkeller and I.\ Pavlyukevich (2009), Limit theorems for $p$-variations of
solutions of SDEs driven by additive stable L\'evy noise and model
selection for paleo-climatic data, In ``Recent Development in
Stochastic Dynamics and Stochastic Analysis", J. Duan, S. Luo and C.
Wang (Eds.), \emph{Interdisciplinary Math. Sci.} \textbf{8}.


\bibitem{Herrchen}M. P. Herrchen (2001),
 \newblock {Stochastic Modeling of Dispersive Diffusion by Non-Gaussian
Noise}, Doctorial Thesis, Swiss Federal Inst. of Tech., Zurich.

\bibitem{Heyde97}C. C.  Heyde (1997),
 \newblock {Quasi-Likelihood and its
Application: A General Approach to Optimal Parameter Estimation}.
Springer, New York.

\bibitem{Horst} W. Horsthemke and R. Lefever (1984),
Noise-Induced Transitions, Springer-Verlag, Berlin.

\bibitem{Hutton1986} J. E. Hutton and P. I. Nelson (1986),
Quasi-likelihood estimation for semimartingales, \emph{Stochastic
Processes and their Applications} \textbf{22}, 245-257.

\bibitem{Ibragimov}I. A. Ibragimov, R. Z. Has'minskii (1981),
\newblock{\em Statistical Estimation-Asymptotic Theory}.
Springer-Verlag.

\bibitem{Ikeda} N. Ikeda and S. Watanabe (1989),
Stochastic Differential Equations and Diffusion Processes,
North-Holland Publishing Company, Amsterdam.

\bibitem{Imkeller} P. Imkeller and I. Pavlyukevich (2002),
Model reduction and stochastic resonance, \textbf{Stochastics and
Dynamics} {\bf 2}(4), 463--506.

\bibitem{ImkellerP-06} P. Imkeller and I. Pavlyukevich (2006),
First exit time of SDEs driven by stable L\'evy processes,
\emph{Stoch. Proc. Appl.} \textbf{116}, 611-642.

\bibitem{ImkellerP-08}
P. Imkeller, I. Pavlyukevich and T. Wetzel (2009), First exit times
for L\'evy-driven diffusions with exponentially light jumps,
\emph{Annals of Probability} \textbf{37}(2), 530¨C564.

\bibitem{Nicolau2004}J. Nicolau (2004),
\newblock{Introduction to the Estimation of Stochastic Differential Equations Based on
Discrete Observations}, Stochastic Finance 2004 (Autumn School and
International Conference).

\bibitem{jacod2006}J. Jacod (2006),
\newblock{Parametric inference for discretely observed non-ergodic
diffusions}, \emph{Bernoulli} \textbf{12}(3), 383-401.

 \bibitem{Aleksander} A. Janicki and A. Weron (1994),
 {Simulation and Chaotic Behavior of $\alpha-$Stable Stochastic
 Processes}, Marcel Dekker, Inc..

\bibitem{Kantz} W. Just, H. Kantz, C. Rodenbeck and M. Helm (2001),
{Stochastic modeling: replacing fast degrees of freedom by noise},
 \emph{J. Phys. A: Math. Gen.} {\bf 34}, 3199--3213.

\bibitem{s13} I. Karatzas and S. E. Shreve (1991), {Brownian Motion and Stochastic
Calculus} 2nd edition, Springer.

\bibitem{Kessler2000} M. Kessler (2000),
{Simple and explicit estimating functions for a discretely observed
diffusion process},
 \emph{Scandinavian Journal of Statistics} \textbf{27}(1), 65-82.

\bibitem{Krishnan} V. Krishnan (2005),
Nonlinear Filtering and Smoothing: An Introduction to Martingales,
Stochastic Integrals and Estimation, Dover Publications, Inc., New
York.


\bibitem{Kleptsyna2000}M. L. Kleptsyna, A. Le Breton and M. C.
Roubaud (2000),
 \newblock {Parameter estimation and optimal filtering for fractional type stochastic systems}.
 \newblock \emph{Statist. Inf. Stochast. Proces.} \textbf{3}, 173-182.

\bibitem{Kleptsyna2002}M. L. Kleptsyna and A. Le Breton (2002),
 \newblock {Statistical analysis of the fractional Ornstein-Uhlenbeck type
 process},
 \newblock \emph{Statistical Inference for Stochastic Processes} \textbf{5}(3), 229-242.


\bibitem{Klebaner}F. Klebaner (2005),
 \newblock {Introduction to Stochastic Calculus with Application},
 \newblock Imperial College Press, Second Edition, 2005.

\bibitem{Kogon}S. Kogon, D. Williams (1998),
\newblock{Characteristic function based estimation of stable distribution
parameters}, \emph{in A practical guide to heavy tails,} M. T. R.
Adler R. Feldman, Ed. Berlin: Birkhauser, 311-335.


\bibitem{Kolmogorov1940}A.N. Kolmogorov (1940),
 \newblock {Wienersche spiralen und einige andere interessante kurven im hilbertschen
 raum},
 \newblock \emph{C.R.(Doklady) Acad. URSS (N.S)} \textbf{26}, 115-118, 1940.


\bibitem{Koutrouvelis}I. A. Koutrouvelis (1980),
\newblock{Regression-type estimation of the parameters of stable laws.}
\emph{Journal of the American Statistical Association} \textbf{75},
918-928.


\bibitem{Kunita2004} H. Kunita (2004), Stochastic differential equations based on L\'evy
processes and stochastic flows of diffeomorphisms, \emph{Real and
stochastic analysis} (Eds. M. M. Rao), 305--373, Birkhuser, Boston,
MA.

\bibitem{Kuruoglu}E. E. Kuruoglu (2001),
\newblock{Density parameter estimationof skewed $\alpha$-stable
distributions}, \emph{Singnal Processing}, IEEE Transactions on
2001, \textbf{49}(10): 2192-2201.

\bibitem{Kutoyants1084b}Yu. A. Kutoyants (1984),
\newblock{Parameter estimation for diffusion type processes of observations},
 \emph{Statistics} \textbf{15}(4), 541-551.

\bibitem{Le Breton}A. Le Breton (1998),
 \newblock {Filtering and parameter estimation in a simple linear model driven by a fractional Brownian
 motion},
 \newblock \emph{Stat. Probab. Lett.} \textbf{38}(3), 263-274.

 \bibitem{Le Breton 1976}A. Le Breton (1976),
 \newblock {On continuous and discrete sampling for parameter estimation in diffusion type
 processes},
 \newblock \emph{Mathematical Programming Study} \textbf{5}, 124-144.

\bibitem{Lipster1977}R. S. Lipster and A. N. Shiryaev (1977),
 \newblock { Statistics of Random Processes},
 \newblock Springer, New York, 1977.


\bibitem{LiuDuan2} X. Liu, J. Duan, J. Liu and  P. E. Kloeden (2009),
 Synchronization of systems of Marcus canonical
equations driven by $\alpha$-stable noises,
 \emph{Nonlinear Analysis - Real World Applications}, to appear, 2009.

\bibitem{A. W. Lo.} A. W. Lo (1991),
 \newblock {Long-term memory in stock market prices},
 \newblock \emph{Econometrica} \textbf{59}, 1279-1313.

\bibitem{B. B. Mandelbrot 1969}B. B. Mandelbrot and J. R. Wallis (1969),
 \newblock {Computer experiments with fractional Gaussian
 noises},
 \newblock \emph{Water Resources Research} \textbf{5}, 228-267.

\bibitem{Mandelbrot and Van Ness} B. B. Mandelbrot and J. W. Van
Ness (1968),
 \newblock {Fractional Brownian motions, fractional noises and
 applications},
 \newblock \emph{SIAM Rev.} \textbf{10}, 422-437.

\bibitem{Mao1995} X. Mao (1995),
 \newblock {Stochastic Differential Equations and Applications},
 Horwood Publishing, Chichester.

 \bibitem{Maslowski} B. Maslowski and B. Schmalfuss (2005),
 Random dynamical systems and stationary solutions of differential
 equationsdriven by the fractional Brownian motion, \emph{Stoch. Anal.
 Appl.} \textbf{22}(6), 1577 - 1607.

\bibitem{McCulloch}J. H. McCulloch (1986),
 \newblock {Simple consistent estimators of stable distributions},
 \newblock \emph{Communications in Statistics-Simulation and Computation} \textbf{15},
 1109-1136.

\bibitem{Yuliya S. Mishura}Y. S. Mishura (2008),
 \newblock {Stochastic Calculus for Fractional Brownian Motion and Related
 Processes},
 \newblock Springer, Berlin.

\bibitem{Moss} F. Moss and P. V. E. McClintock (eds.),
\textit{Noise in Nonlinear Dynamical Systems}. Volume 1: Theory of Continuous Fokker-Planck Systems (2007);
Volume 2: Theory of Noise Induced Processes in Special Applications (2009); Volume 3: Experiments and Simulations (2009).
   Cambridge University Press.

\bibitem{Nolan} J. P. Nolan (2007),
\emph{Stable Distributions - Models for Heavy Tailed Data},
Birkh\"ause, Boston, 2007.

\bibitem{Nourdin}
I. Nourdin and T. Simon (2006), \newblock {On the absolute
continuity of L\'{e}vy processes with drift}, \emph{Ann. Prob.}
\textbf{34}(3), 1035-1051.

\bibitem{Norros1999}I. Norros, E. Valkeila and J. Virtamo (1999),
 \newblock {An elementary approach to a Girsanov formula and other analytical results on fractional Brownian
 motions},
 \newblock \emph{Bernoulli} \textbf{5}(4), 571-587.

 \bibitem{Nualart} D. Nualart (2003), Stochastic calculus with respect to the fractional
Brownian motion and applications, \emph{Contemporary Mathematics},
\textbf{336}, 3-39.


\bibitem{Oksendal3} B. Oksendal (2005),
Applied Stochastic Control Of Jump Diffusions, Springer-Verlag, New
York.

\bibitem{Oksendal}B. Oksendal (2003),
\newblock{Stochastic Differenntial Equations}, Sixth Ed.,
  Springer-Verlag, New York.

\bibitem{Bernt Oksendal} B. Oksendal, F. Biagini, T. Zhang and Y. Hu
(2008),
 \newblock {Stochastic Calculus for Fractional Brownian Motion and Applications}.
 \newblock Springer.

\bibitem{Palmer2} T. N. Palmer,  G. J. Shutts, R. Hagedorn,
F. J. Doblas-Reyes, T. Jung and M. Leutbecher (2005), Representing
model uncertainty in weather and climate prediction, \emph{Annu.
Rev. Earth Planet. Sci.} \textbf{33}, 163-193.

\bibitem{Papoulis} A. Papoulis (1984),
Probability, Random Variables, and Stochastic Processes, McGraw-Hill
Companies, 2nd edition.

\bibitem{Pearson1994}N. D. Pearson and T. Sun (1994),
 \newblock {Exploiting the conditional density in estimating the term structure: an application to the Cox, Ingersoll and Ross
 model},
 \newblock \emph{The Journal of Finance} \textbf{49}(4), 1279-1304.

\bibitem{Pedersen}A. R. Pedersen (1995),
Consistency and asymptotic normality of an approximate maximum
likelihood estimator for discretely observed diffusion
 processes,
 \newblock \emph{Bernoulli} \textbf{1}(3), 257-279.

\bibitem{C. K. Peng}C. K. Peng, V. Buldyrev, S. Havlin, M. Simons, H. E. Stanley, and
A. L. Goldberger (1994),
 \newblock {Mosaic organization of DNA nucleotides}.
 \newblock \emph{Phys. Rev. E} \textbf{49}, 1685-1689.

\bibitem{PZ} S. Peszat and J. Zabczyk (2007),
Stochastic Partial Differential Equations with L{\'e}vy Processes,
Cambridge University Press, Cambridge, UK.

\bibitem{prakasa1999a}
B.L.S. Prakasa Rao (1999),
\newblock {Statistical Inference for Diffusion Type Processes},
Arnold, London.

\bibitem{prakasa}
B.L.S. Prakasa Rao (1999),
\newblock {Semimartingales and their Statistical Inference},
\newblock {Chapman \& Hall/CRC}.

\bibitem{B. L. S. Prakasa Rao2003}B. L. S. Prakasa Rao (2003),
 \newblock {Parametric estimation for linear stochastic differential equations driven by fractional Brownian motion}.
 \newblock http://www.isid.ac.in/~statmath/eprints


\bibitem{Press 1972}S. Press (1972),
\newblock {Estimation of univariate and multivariate stable distributions},
\emph{Journal of the Americal Statistical Association} \textbf{67},
842-846.

\bibitem{Protter} P. E. Protter (2005), Stochastic Integration and Differential
Equations, Springer-Verlag, New York, Second Edition.

\bibitem{Roz} B. L. Rozovskii (1990),
Stochastic Evolution Equations, Kluwer Academic Publishers, Boston.

\bibitem{Robinson1977} P. M. Robinson (1977),
Estimation of a time series model from unequally spaced data,
\emph{Stoch. Proc. Appl.} \textbf{6}, 9-24.

\bibitem{Gennady Samorodnitsky2000}G. Samorodnitsky, M. S. Taqqu (2008),
 \newblock {\em Stable Non-Gaussian Random Processes- Stochastic Models with Infinite Variance}.
 \newblock Chapman \& Hall/CRC.

\bibitem{Sato-99}K. Sato (1999),
 \newblock {L\'evy  Processes and Infinitely Divisible Distrributions,}
Cambridge University Press, Cambridge, UK, 1999

\bibitem{SSB91}H. Scher, M. F. Shlesinger and J. T. Bendler (1991),
 \newblock{Time-scale invariance in transport and relaxation}, \emph{Phys.
 Today} \textbf{44}(1), 26-34.

 \bibitem{Schertzer}
D. Schertzer, M. Larcheveque, J. Duan,  V. Yanovsky and S. Lovejoy
(2000),
\newblock{Fractional Fokker--Planck equation for non-linear
stochastic differential equations driven by non-Gaussian L\'evy stable
noises}, {\em J. Math. Phys.} {\bf 42}, 200-212.

\bibitem{Shlesinger}M. F. Shlesinger, G. M. Zaslavsky and U. Frisch (1995),
 \newblock{L\'evy Flights and Related Topics in Physics},
Lecture Notes in Physics, Springer-Verlag, Berlin.

\bibitem{Sorensen1998b} M. Sorensen (1999),
\newblock {On asymptotics of estimating functions},
\emph{Brazillian Journal of Probability and Statistics} \textbf{13},
111-136.

\bibitem{Stroock1979}D. W. Stroock and S. R. S. Varadhan (1979),
 \newblock {Multidimensional Diffusion Processes},
 Springer Verlag, Berlin.

\bibitem{Swinney} T. H. Solomon, E. R. Weeks, and H. L. Swinney (1993),
Observation of anomalous diffusion and L\'evy flights in a
two-dimensional rotating flow, \emph{Phys. Rev. Lett.}
\textbf{71}(24), 3975 - 3978.

\bibitem{M. S. Taqqu 1997}M. S. Taqqu, V. Teverovsky, and W.
Willinger (1995),
 \newblock {Estimators for long-range dependence: an empirical
 study},
 \newblock \emph{Fractals}, \textbf{3}(4), 785-798.

\bibitem{Tsihrintzis1996}G. A. Tsihrintzis and C.
L. Nikias (1995),
 \newblock {Fast estimation of the parameters of alpha-stable impulsive
 interference using asymptotic extreme value theory},
 \newblock ICASSP-95, \textbf{3},
 1840-1843.


 \bibitem{VanKampen2} N. G. Van Kampen (1987),
How do stochastic processes enter into physics? Lecture Note in
Mathe. {\bf 1250/1987}, 128--137.

\bibitem{VanKampen3} N. G. Van Kampen (1981), Stochastic Processes in Physics and
Chemistry, North-Holland, New York.

\bibitem{WaymireDuan} E. Waymire  and J. Duan (Eds.) (2005),
Probability and Partial Differential Equations in Modern Applied
Mathematics, Springer-Verlag.

\bibitem{Wilks} D. S. Wilks (2005), Effects of stochastic
parameterizations in the Lorenz '96 system, \emph{Q. J. R. Meteorol.
Soc.} \textbf{131}, 389-407.

\bibitem{Williams}P. D. Williams (2005), Modeling climate change: the
role of unresolved processes, \emph{Phil. Trans. R. Soc. A}
\textbf{363}, 2931-2946.

\bibitem{Wong} E. Wong and B. Hajek (1985), Stochastic Processes in
Engineering Systems, Spring-Verlag, New York.

\bibitem{Woy} W. A. Woyczynski (2001), L\'evy processes in the physical
sciences, In L\'evy processes: theory and applications, O. E.
Barndorff-Nielsen, T. Mikosch and S. I. Resnick (Eds.), 241-266,
Birkh\"auser, Boston, 2001.


\bibitem{Yaglom1958} A. M. Yaglom (1958),
 \newblock {Correlation theory of processes with random stationary nth
 increments},
 \newblock \emph{AMS Transl.} \textbf{2}(8), 87-141.

\bibitem{YangDuan}
Z. Yang and J. Duan (2008), An intermediate regime for exit
phenomena driven by non-Gaussian L\'evy noises, \emph{Stochastics
and Dynamics} \textbf{8}(3), 583-591.

\bibitem{Yon96} F. Yonezawa (1996), Introduction to focused session on
`anomalous relaxation, \emph{J. Non-Cryst. Solids} \textbf{198-200},
503-506.

\bibitem{Yoshida} N. Yoshida (2004), {\em Estimation for diffusion processes from discrete observations},
\emph{J. Multivariate Anal.} \textbf{41}(2), 220-242.

\end{thebibliography}
\end{document}